
\input amstex
\documentstyle{amsppt}
\magnification=1200
\catcode`\@=11
\redefine\logo@{}
\catcode`\@=13
\pageheight{19cm}

\define \bn{\Bbb N}
\define \bz{\Bbb Z}
\define \bq{\Bbb Q}
\define \br{\Bbb R}
\define \bc{\Bbb C}

\define \bp{\Bbb P}
\define \bk{\Bbb K}

\define\Da{{\Cal D}}

\define\La{{\Cal L}}
\define\Ka{{\Cal K}}

\define\rk{\text{rk}~}


\define\pr{\text{pr}}

\define\Exc{\text{Exc}}
\define\0o{{\overline 0}}
\define\1o{{\overline 1}}

\TagsOnRight

\document

\topmatter

\title
On a classical correspondence between K3 surfaces
\endtitle

\author
Carlo Madonna and
Viacheslav V. Nikulin \footnote{Supported by Russian Fund of
Fundamental Research (grant N 00-01-00170).
\hfill\hfill}
\endauthor

\address
Dottorato di Ricerca, Dipartimento di Matematica, Univ. Roma
2 Tor Vergata, Roma
\endaddress
\email
madonna\@mat.uniroma3.it
\endemail

\address
Deptm. of Pure Mathem. The University of Liverpool, Liverpool
L69 3BX, UK;
\vskip1pt
Steklov Mathematical Institute,
ul. Gubkina 8, Moscow 117966, GSP-1, Russia
\endaddress
\email
vnikulin\@liv.ac.uk\ \
slava\@nikulin.mian.su
\endemail

\dedicatory
To 80-th Birthday of I.R. Shafarevich
\enddedicatory

\abstract
Let $X$ be a K3 surface which is intersection of three
(i.e. a net ${\Bbb P}^2$) of quadrics in ${\Bbb P}^5$.
The curve of degenerate quadrics has degree 6 and
defines a natural double covering $Y$ of ${\Bbb P}^2$ ramified
in this curve which is again a K3.
This is a classical example of a correspondence between
K3 surfaces which is related with moduli of sheaves on K3's studied
by Mukai. When general (for fixed Picard lattices) $X$ and $Y$ are
isomorphic? We give necessary and sufficient conditions in terms of
Picard lattices of $X$ and $Y$.

E.g. for Picard number 2 the Picard lattice of $X$ and $Y$ is defined by
its determinant $-d$ where $d>0$, $d\equiv 1 \mod 8$, and one of equations
$a^2-db^2=8$ or $a^2-db^2=-8$ has an integral solution $(a,b)$.
Clearly, the set of these $d$ is infinite: $d\in \{(a^2\mp 8)/b^2\}$
where $a$ and $b$ are odd integers. This gives all possible divisorial
conditions on the 19-dimensional moduli of intersections of three quadrics
$X$ in ${\Bbb P}^5$ which imply $Y\cong X$. 
One of them, when $X$ has a line is classical and corresponds to $d=17$.

Similar considerations can be applied to a realization
of an isomorphism
$(T(X)\otimes \bq, H^{2,0}(X)) \cong (T(Y)\otimes \bq, H^{2,0}(Y))$
of transcendental periods over $\bq$ of two K3 surfaces $X$ and $Y$
by a fixed sequence of types of Mukai vectors.
\endabstract

\rightheadtext
{Correspondences for K3}
\leftheadtext{C. Madonna and V.V. Nikulin}
\endtopmatter

\document

\head
0. Introduction
\endhead

In this paper we study a classical correspondence between algebraic
K3 surfaces over $\bc$.

Let a K3 surface $X$ be an intersection of three quadrics in ${\Bbb P}^5$
(more generally, $X$ is a K3 surface with a primitive
polarization $H$ of degree 8).
The three quadrics define the projective plane ${\Bbb P}^2$ (the net)
of quadrics. Let $C\subset {\Bbb P}^2$ be the curve of degenerate
quadrics. The curve $C$ has degree 6 and defines another
K3 surface $Y$  which is the minimal resolution of singularities of
the double covering of ${\Bbb P}^2$ ramified in $C$. It has the natural
linear system $|h|$ with $h^2=2$ which is preimage of lines in $\bp^2$.
This is a classical and a very beautiful example of a correspondence between
K3 surfaces. It is defined by a 2-dimensional algebraic cycle
$Z\subset X\times Y$.

This example is related with the moduli of sheaves on K3 surfaces
studied by Mukai \cite{5}, \cite{6}. It is well-known that the K3
surface $Y$ is the moduli of sheaves ${\Cal E}$ on $X$ with the rank
$r=2$,
the first Chern class $c_1({\Cal E})=H$ and the Euler characteristic
$\chi=\chi({\Cal E})=4$.
We apply this construction to study the following questions.

\proclaim{Question 1} When $Y$ is isomorphic to $X$?
\endproclaim

We want to answer this question in terms of the Picard lattices
$N(X)$ and $N(Y)$ of
$X$ and $Y$. Then our question is as follows:

\proclaim{Question 2} Assume that $N$ is a hyperbolic lattice,
$\widetilde{H}\in N$ a primitive element with square $8$. What are
conditions on $N$ and $\widetilde{H}$ such that for any K3 surface
$X$ with the Picard lattice $N(X)$ and a primitive polarization $H\in N(X)$
of degree 8 the corresponding K3 surface $Y$ is isomorphic to $X$
if the the pairs of lattices $(N(X), H)$ and
$(N, \widetilde{H})$ are isomorphic as abstract lattices with
fixed elements?

In other words, what are conditions on $(N(X), H)$ as an abstract lattice
with a primitive vector $H$ with $H^2=8$ which are sufficient for 
$Y$ to be isomorphic to $X$ and are necessary  
if $X$ is a general K3 surfaces with the Picard lattice $N(X)$?

\endproclaim

We give an answer to the questions 2 in Theorems 2.2.3 and 2.2.4 and
also Propositions 2.2.1 and 2.2.2. In particular, if the Picard number
$\rho (X)=\rk N(X) \ge 12$, the result is very simple:
$X\cong Y$, if and only if there exists
$x\in N(X)$ such that $x\cdot H\equiv 1\mod 2$.
This follows from results of Mukai \cite{6} and also \cite{7},
\cite{8}.

K3 surfaces $X$ and $Y$ with
$\rho (X)=\rho (Y)=2$ are especially interesting.
Really, it is well-known that the moduli space of K3 surfaces
$X$, which are intersections of three quadrics, is $19$-dimensional.
If $X$ is general, i.e. $\rho(X)=1$, then the surface
$Y$ cannot be isomorphic to $X$
because $N(X)=\bz H$ where $H^2=8$ and $N(X)$ does not have elements
with square $2$ which is necessary if $Y\cong X$. Thus, if $Y\cong X$,
then $\rho (X)\ge 2$, and $X$ belongs to a codimension 1 submoduli
space of K3 surfaces which is a
{\it divisor in the 19-dimensional moduli space of
intersections of three quadrics in $\bp^5$ (up to codimension 2).}
To describe connected components of this divisor, it is equivalent
to describe Picard lattices $N(X)\cong N(Y)$ of the surfaces
$X\cong Y$ above with $\rho (X)=\rho (Y)=2$ such that general K3 surfaces
$X$ and $Y$ with these Picard lattices have $X\cong Y$. We show that
the Picard lattice $N(X)\cong N(Y)$ of these $X\cong Y$ is defined by its
determinant  which is equal to $-d$ where $d>0$ and $d\equiv 1\mod 8$.
We show that the set $\Da$ of these numbers $d$ is exactly the set of
$d\in \bn$ such that $d\equiv 1\mod 8$ and one of the equations
$$
a^2-db^2=8
\tag{0.1}
$$
or
$$
a^2-db^2=-8
\tag{0.2}
$$
has an integral solution. It is easy to see that solutions $(a,b)$ of
these equations are odd. It follows that
$\Da=\Da_+\cup \Da_-$ where
$$
\Da_+ =\{{a^2-8\over b^2}\in \bn\ |\ a,\ b\in \bn\
\text{are odd}\}
\tag{0.3}
$$
and
$$
\Da_- =
\{{a^2+8\over b^2}\in \bn\ |\ a,\ b\in \bn\ \text{are odd}\}.
\tag{0.4}
$$
Both sets $\Da_+$ and $\Da_-$ are infinite. The set $\Da_+$ contains
the infinite sequence $a^2-8$ where $a\in \bn$
is odd. The set $\Da_-$ contains the infinite sequence $a^2+8$
where $a\in \bn$ is odd. Similarly, we show that the sets
$\Da_+-\Da_+\cap \Da_-$ and $\Da_- - \Da_+\cap \Da_-$ are infinite
(see considerations after Theorem 3.1.7).
In contrary, we don't know if the set
$\Da_+\cap \Da_-$ is infinite.
In Sect. 3.2, we give an algorithm and a program for
calculation of $\Da$ and describe geometry of the surfaces $Y\cong X$
for $d\in \Da$. See Theorems 3.2.1 --- 3.2.5 and also results of
calculations (for small $d$) after Theorem 3.2.5. Calculations using
this algorithm give the list of first numbers
from $\Da$: 1, 9, 17, 33, 41, 57, 73, 89, 97, 113, ..., 2009. 
See Theorem 3.2.5.

The set $\Da$ labels connected components
of the divisor, where $Y\cong X$,
in 19-dimensional moduli of intersections of
three quadrics $X$. Each $d\in \Da$, gives a connected 18-dimensional
moduli space of K3 surfaces with the Picard lattice
$N(X)=N(Y)$ of rank 2 and determinant $-d$ where $X\cong Y$.
E.g. it is well-known that $Y\cong X$ if $X$ has a line.
This is a divisorial condition on moduli of $X$.
This component is labeled by $d=17\in \Da$.
Thus, our results show that: {\it ``There exists
an infinite number of different divisorial conditions on moduli of
intersections of three quadrics $X$ in $\bp^5$ such that each of them 
implies $Y\cong X$. They are labeled by elements of the infinite 
set $\Da$ which  was described above. The case $d=17$ corresponds 
to the classical example above of $X$ containing a line.''}

We mention that solutions $(a,b)$ of the equations
\thetag{0.1} and \thetag{0.2} can be interpreted as elements of
Picard lattices of $X$ and $Y$. E. g. for general
$X$ and $Y$ with $\rho (X)=\rho(Y)=2$,
we have $Y\cong X$ if and only if $\det N(X)$ is odd
and there exists $h_1\in N(X)$ such that $(h_1)^2=\pm 4$ (for one
of signs) and
$h_1\cdot H \equiv 0\mod 2$. Similar condition in terms of $Y$ is:
$\det N(X)$ is odd and there exists $h_1\in N(Y)$ such that
$h_1^2=\pm 4$ (for one of signs). It follows the following very
simple sufficient condition when $Y\cong X$ (see Corollary 3.1.9)
which we want to formulate exactly.

\proclaim{Theorem}
Let $X$ be a K3 surface.
Assume that $X$ is an intersection of three quadrics (more generally,
$X$ has a primitive polarization $H$ of degree 8).
Let $Y$ be a K3 surface which is the double covering
of the net $\bp^2$ of the quadrics defining $X$  
ramified along the curve of degenerate quadrics 
(more generally, $Y$ is the moduli space of sheaves
on $X$ with the Mukai vector $v=(2,H,2)$). 

Then $Y\cong X$, if there exists $h_1\in N(X)$ such that
the primitive sublattice $[H,h_1]_{\text{\pr}}$ in $N(X)$
generated by $H$ and $h_1$ has odd determinant, and
$$
(h_1)^2=\pm 4\ \text{and\ } h_1\cdot H\equiv 0\mod 2.
$$
These conditions are necessary if either $\rho(X)=1$, or $\rho(X)=2$ and
$X$ is a general K3 surface with its Picard lattice.
\endproclaim

From our point of view, this statement is very interesting
because elements $h_1$ of the Picard lattice $N(X)$ with negative
square $(h_1)^2=-4$ get some geometrical meaning. For K3 surfaces
it is well-known only for elements $\delta$ of the Picard lattice
$N(X)$ with negative square $\delta^2=-2$: then
$\delta$ or $-\delta$ is effective.

It seems, many known examples of $Y\cong X$ (e. g. see \cite{3}, 
\cite{15}) follow from the theorem. 

\smallpagebreak 

Similar methods can be developed for much more general situation.
Let $X$ and $Y$ are K3 surfaces,
$$
\phi:(T(X)\otimes \bq, H^{2,0}(X))\cong (T(Y)\otimes \bq, H^{2,0}(Y))
\tag{0.5}
$$
an isomorphism of their transcendental periods over $\bq$, and
$$
(a_1,H_1,b_1)^{\pm}, \dots , (a_k,H_k,b_k)^{\pm},
\tag{0.6}
$$
a sequence of types of isotropic Mukai vectors of moduli of sheaves
on K3 where $\pm$ shows the direction of the correspondence.

Similar methods and calculations can be applied to study the
following question:
{\it When there exists a correspondence between $X$ and
$Y$ which is given by the sequence \thetag{0.6} of Mukai vectors, and
which gives the isomorphism \thetag{0.5} between their transcendental
periods?}

In \cite{6} and \cite{9} sufficient and necessary  
conditions on \thetag{0.5} were given when there exists at 
least one such a sequence \thetag{0.6} 
with coprime Mukai vectors $(a_i,H_i,b_i)$.

\smallpagebreak

The fundamental tool to get the results above is the Global
Torelli Theorem for K3 surfaces proved by I.I. Piatetskii-Shapiro
and I.R. Shafarevich in \cite{10}: Using results of Mukai
\cite{5}, \cite{6}, we can
calculate periods of $Y$ using periods of $X$; by the Global Torelli
Theorem \cite{10}, we can find out if $Y$ is isomorphic to $X$.

\smallpagebreak

We are grateful to A.N. Tyurin and A. Verra for useful and stimulating
discussions.

\smallpagebreak

\head
1. Preliminary notations and results about lattices and K3 surfaces
\endhead

\subhead
1.1. Some notations about lattices
\endsubhead
We use notations and terminology from \cite{8} about lattices,
their discriminant groups and forms. A {\it lattice} $L$ is a
non-degenerate integral symmetric bilinear form. I. e. $L$ is a free
$\bz$-module equipped with a symmetric pairing $x\cdot y\in \bz$ for
$x,\,y\in L$, and this pairing should be non-degenerate. We denote
$x^2=x\cdot x$. The {\it signature} of $L$ is the signature of the
corresponding real form $L\otimes \br$.
The lattice $L$ is called {\it even}
if $x^2$ is even for any $x\in L$. Otherwise, $L$ is called {\it odd}.
The {\it determinant} of $L$ is defined to be $\det L=\det(e_i\cdot e_j)$
where $\{e_i\}$ is some basis of $L$. The lattice $L$ is {\it unimodular}
if $\det L=\pm 1$.

The {\it dual lattice} of $L$ is
$L^\ast=Hom(L,\,\bz)\subset L\otimes \bq$. The
{\it discriminant group} of $L$ is $A_L=L^\ast/L$. It has the order
$|\det L|$. The group $A_L$ is equipped with the
{\it discriminant bilinear form} $b_L:A_L\times A_L\to \bq/\bz$
and the {\it discriminant quadratic form} $q_L:A_L\to \bq/2\bz$
if $L$ is even. To get this forms, one should extend the form of $L$ to
the form on the dual lattice $L^\ast$ with values in $\bq$.

For $x\in L$, we shall consider
the invariant $\gamma(x)\ge 0$ where
$$
x\cdot L=\gamma (x)\bz .
\tag{1.1.1}
$$
Clearly, $\gamma (x)|x^2$ if $x\not=0$. 

We denote by $L(k)$ the lattice obtained from a lattice $L$ by
multiplication of the form of $L$ by $k\in \bq$.

The orthogonal sum of lattices $L_1$ and $L_2$ is denoted by
$L_1\oplus L_2$.

For a symmetric integral matrix
$A$, we denote by $\langle A \rangle$ a lattice which is given by
the matrix $A$ in some bases. We denote
$$
U=\left\langle
\matrix
0&1\\
1&0
\endmatrix
\right\rangle.
\tag{1.1.2}
$$
Any even unimodular lattice of the signature $(1,1)$ is isomorphic to
$U$.

An embedding $L_1\subset L_2$ of lattices is called {\it primitive}
if $L_2/L_1$ has no torsion.

We denote by $O(L)$, $O(b_L)$ and $O(q_L)$ the automorphism groups of
the corresponding forms. Any $\delta\in L$ with $\delta^2=-2$ defines
a reflection $s_\delta\in O(L)$ which is given by the formula
$x\to x+(x\cdot \delta)\delta$, $x\in L$. All such reflections generate
the {\it 2-reflection group} $W^{(-2)}(L)\subset O(L)$.

\subhead
1.2. Some notations about K3 surfaces
\endsubhead
Here we remind some basic notions and results about K3 surfaces,
e. g. see \cite{10}, \cite{11}, \cite{12}.
A K3 surface $S$ is a projective algebraic surface over
$\bc$ such that its canonical class $K_S$ is zero and the irregularity
$q_S=0$. We denote by $N(S)$ the {\it Picard lattice} of $S$ which is
a hyperbolic lattice with the intersection pairing
$x\cdot y$ for $x,\,y\in N(S)$. Since the canonical class $K_S=0$,
the space $H^{2,0}(S)$ of 2-dimensional holomorphic differential
forms on $S$ has dimension one over $\bc$, and
$$
N(S)=\{x\in H^2(S,\bz)\ |\ x\cdot H^{2,0}(S)=0\}
\tag{1.2.1}
$$
where $H^2(S,\bz)$ with the intersection pairing is a
22-dimensional even unimodular lattice of signature
$(3,19)$. The orthogonal lattice $T(S)$ to $N(S)$ in $H^2(S,\bz)$ is called
the {\it transcendental lattice of $S$.} We have
$H^{2,0}(S)\subset T(S)\otimes \bc$. The pair $(T(S), H^{2,0}(S))$ is
called the {\it transcendental periods of $S$}.
The {\it Picard number} of $S$ is
$\rho(S)=\rk N(S)$. An non-zero element $x\in N(S)\otimes \br$ is
called {\it nef} if $x\not=0$ and $x\cdot C\ge 0$ for any effective
curve $C\subset S$. It is known that an element $x\in N(S)$ is ample
if $x^2>0$, $x$ is $nef$, and the orthogonal complement
$x^\perp$ to $x$ in  $N(S)$ has no elements with square $-2$.
For any element $x\in N(S)$ with $x^2\ge 0$, there exists a reflection
$w\in W^{(-2)}(N(S))$ such that the element
$\pm w(x)$ is nef; it then is ample, if $x^2>0$ and $x^\perp$ had no
elements with square $-2$ in $N(S)$.

We denote by $V^+(S)$ the light cone of $S$, which is the half-cone of
$$
V(S)=\{x\in N(S)\otimes \br\ |\ x^2>0\ \}
\tag{1.2.2}
$$
containing a polarization of $X$. In particular, all $nef$ elements
$x$ of $S$ belong to $\overline{V^+(S)}$:
one has $x\cdot V^+(S)>0$ for them.

The reflection group $W^{(-2)}(N(S))$ acts in $V^+(S)$ discretely,
and its  fundamental
chamber is the closure $\overline{\Ka(S)}$ of the K\"ahler cone
$\Ka(S)$ of $S$. It is the same as the set of all $nef$ elements of $S$.
Its faces are orthogonal to the set $\Exc(S)$ of all exceptional curves $r$
on $S$ which are non-singular rational curves $r$ on $S$ with $r^2=-2$.
Thus, we have
$$
\overline{\Ka (S)}=\{0\not=x\in \overline{V^+(S)}\ |
\ x\cdot \Exc(S)\ge 0\,\}.
\tag{1.2.3}
$$

\subhead
1.3. K3 surfaces with polarizations of degree 8 and 2
\endsubhead
The following results are well-known (see Mayer \cite{4},
Saint-Donat \cite{11} and Shokurov \cite{13}).

Let $X$ be a K3 surface with a primitive polarization
$H\in N(X)$ of degree $H^2=8$. Here primitive means that the sublattice
$\bz H\subset N(X)$ is primitive.

\proclaim{Proposition 1.3.1} The linear system $|H|$ has dimension $5$,
and there are the following and the only following
cases for the linear system $|H|$:

(a) $|H\cdot E|>2$ for any elliptic curve $E$ on $X$ (it has $E^2=0$).
Then the linear system $|H|$ gives an embedding of $X$ to $\bp^5$ as
intersection of three quadrics.

(b) $|H\cdot E|\ge 2$ for any elliptic curve $E$ on $X$, and there
exists an elliptic curve $E$ on $X$ such that $H\cdot E=2$. Then the
linear system $|H|$ is hyperelliptic and gives a double covering
of a rational scroll. In this case, $H$ and $E$ generate a
primitive sublattice of $N(X)$ which is isomorphic to $U(2)$. 

(c) $|H\cdot E|\ge 1$ for any elliptic curve $E$ on $X$, and there
exists an elliptic curve $E$ on $X$ such that $H\cdot E=1$.
Then the linear system $|H|$ has a fixed component $D$ which is
a non-singular rational curve on $X$ (it has $D^2=-2$),
and $|H|=5|E|+D$ where $|E|$
is an elliptic pencil on $X$ and $E\cdot D=1$.
In this case, $H$ and $E$ generate
a primitive sublattice of $N(X)$ which is isomorphic to $U$. 
\endproclaim

Since $H^2=8$, we then have $\gamma (H)=1$, $2$, $4$ or $8$ 
(for $H\in N(X)$). 

For the case (c), $\rho(X)\ge 2$ and $\gamma(H)=1$.

For the case (b), $\rho(X)\ge 2$ and $\gamma(H)=1$ or $2$.
If $\rho(X)=2$, then $\gamma(H)=2$.

For the case (a), $\rho (X)\ge 1$ and $\gamma(H)=1$, $2$, $4$ or $8$.
If $\rho (X)=1$, then $\gamma(H)=8$.

A general K3 surface $X$ has $\rho(X)=1$, then one has the case
(a) and $X$ is an intersection of three quadrics.

\smallpagebreak

Now let $Y$ be a K3 surface with a $nef$ element $h\in N(Y)$ of degree
$h^2=2$. Obviously, $h$ is primitive.

\proclaim{Proposition 1.3.2} The linear system $|h|$ has dimension $2$,
and there are the following and the only following
cases:

(a) $|h\cdot E|\ge 2$ for any elliptic curve $E$ on $Y$.  Then the
linear system $|h|$ gives a double covering of $\bp^2$
ramified along a curve of degree 6 with at most double singularities.

(b) $|h\cdot E|\ge 1$ for any elliptic curve $E$ on $Y$,
and there exists an elliptic curve $E$ on $Y$ such that $h\cdot E=1$.
Then the linear system $|h|$ has a fixed component $D$ which is a
non-singular rational curve (it has $D^2=-2$), and $|h|=2|E|+D$
where $|E|$ is an elliptic pencil and $E\cdot D=1$. In this case,
$h$ and $E$ generate a primitive sublattice of $N(X)$ which is isomorphic to
$U$.
\endproclaim

Since $h^2=2$, we have $\gamma (h)=1$ or $2$.

For the case (b), $\rho(Y)\ge 2$ and $\gamma(h)=1$.

For the case (a), $\rho (Y)\ge 1$ and $\gamma(h)=1$ or $2$.
If $\rho (Y)=1$, then $\gamma(h)=2$.

A general K3 surface $Y$ has $\rho(Y)=1$, then $\gamma(h)=2$,
and one has the case (a).

\head
2. General results on the classical correspondence between K3 surfaces
with primitive polarizations of degree 8 and 2 which gives isomorphic K3's
\endhead

\subhead
2.1. The correspondence
\endsubhead
Let a K3 surface $X$ be an intersection of three quadrics in ${\Bbb P}^5$
(more generally, $X$ is a K3 surface with a primitive
polarization $H$ of degree 8). See general results about intersections
of quadrics in \cite{14}.
The three quadrics define the projective plane ${\Bbb P}^2$ (the net)
of quadrics. Let $C\subset {\Bbb P}^2$ be the curve of
degenerate quadrics. The curve $C$ has degree 6 and defines another
K3 surface $Y$  which is the minimal resolution of singularities of
the double covering of ${\Bbb P}^2$ ramified in $C$. It has the natural
linear system $|h|$ with $h^2=2$ which is preimage of lines on $\bp^2$.
This is a classical and a very beautiful example of a correspondence between
K3 surfaces. It is defined by a 2-dimensional algebraic cycle
$Z\subset X\times Y$.

This example is related with the moduli of sheaves on K3 surfaces
studied by Mukai \cite{5}, \cite{6}. It is well-known that the K3
surface $Y$ is the moduli of sheaves ${\Cal E}$ on $X$ with rank
$r=2$, first Chern class $c_1({\Cal E})=H$ and Euler characteristic
$\chi=\chi({\Cal E})=4$.

Let
$$
H^\ast(X,\bz)=H^0(X,\bz)\oplus H^2(X,\bz)\oplus H^4(X,\bz)
\tag{2.1.1}
$$
be the full cohomology lattice of $X$ with the Mukai product
$$
(u,v)=-(u_0\cdot v_2+u_2\cdot v_0)+u_1\cdot v_1
\tag{2.1.2}
$$
for $u_0,v_0\in H^0(X,\bz)$, $u_1,v_1\in H^2(X,\bz)$,
$u_2,v_2\in H^4(X,\bz)$. We naturally identify
$H^0(X,\bz)$ and $H^4(X,\bz)$ with $\bz$. Then the Mukai product is
$$
(u,v)=-(u_0v_2+u_2v_0)+u_1\cdot v_1.
\tag{2.1.3}
$$

The element
$$
v=(2,H,2)=(r,H,\chi-r)\in H^\ast(X,\bz)
\tag{2.1.4}
$$
is called {\it the Mukai vector} of the net of quadrics or
(more generally) of sheaves on $X$ for which $Y$ is
the moduli space. It is isotropic, i.e. $v^2=0$.
Mukai \cite{5}, \cite{6} showed that one has the natural
identification
$$
H^2(Y,\,\bz)= (v^\perp/\bz v)
\tag{2.1.5}
$$
which also gives the isomorphism of the Hodge structures of $X$ and $Y$.
The element $h=(-1,0,1)\in v^\perp$ has square $h^2=2$, $h\mod \bz v$ belongs
to the Picard lattice $N(Y)$ of $Y$, and the linear system $|h|$
defines the structure of double plane $Y$ for the net of quadrics defining
$X$.

We apply this construction to study the following questions.

\proclaim{Question 2.1.1} When $Y$ is isomorphic to $X$?
\endproclaim

We want to answer this question in terms of Picard lattices of
$X$ and $Y$. Then the exact formulation of our question is
as follows:

\proclaim{Question 2.1.2} Assume that $N$ is a hyperbolic lattice,
$\widetilde{H}\in N$ a primitive element with square $8$. What are
conditions on $N$ and $\widetilde{H}$ such that for any K3 surface
$X$ with the Picard lattice $N(X)$ and a primitive polarization
$H\in N(X)$ of degree 8 the corresponding K3 surface $Y$ is isomorphic
to $X$ if the the pairs of lattices $(N(X), H)$ and
$(N, \widetilde{H})$ are isomorphic as abstract lattices with
fixed elements?

In other words, what are conditions on $(N(X), H)$ as an abstract lattice
with a primitive vector $H$ with $H^2=8$ which are sufficient for 
$Y$ to be isomorphic to $X$ and are necessary if $X$ is
a general K3 surface with the Picard lattice $N(X)$?
\endproclaim

We formulate the results below.

\subhead
2.2. Formulation of general results
\endsubhead
Below $X$ is a K3 surface with a primitive polarization $H$ of degree
$H^2=8$, and $Y$ the corresponding K3 surface $Y$ with the nef element $h$
of degree $h^2=2$ defined in Sect. 2.1.

The following statement follows from the Mukai identification
\thetag{2.1.5} and results from \cite{8}.

\proclaim{Proposition 2.2.1} If $Y$ is isomorphic to $X$, then
the invariant $\gamma (H)$ for $H$ in $N(X)$
(see \thetag{1.1.1}) is equal to $1$.

Assume that $\gamma (H)=1$ for $H$ in $N(X)$.
Then the Mukai identification \thetag{2.1.5}
canonically identifies the transcendental periods
$(T(X), H^{2,0}(X))$ and
\newline
$(T(Y), H^{2,0}(Y))$. It follows that
the Picard lattices $N(Y)$ and $N(X)$ have the same genus.
In particular, $N(Y)$ is isomorphic to $N(X)$ if the genus of $N(X)$
contains only one class.  If the genus of $N(X)$ contains only one class,
then $Y$ is isomorphic to $X$, if additionally the canonical homomorphism
$O(N(X))\to O(q_{N(X)})$ is epimorphic.
Both these conditions are valid (in particular, $Y\cong X$),
if either $\rho(X)\ge 12$ or $N(X)$ contains a copy of $U$
(see \thetag{1.1.2}).
\endproclaim

From now on we can assume that $\gamma(H)=1$ in $N(X)$ since only
for this case we may have $Y\cong X$.

Calculations below are valid for an arbitrary K3 surface $X$ and
a primitive vector $H\in N(X)$ with $H^2=8$ and $\gamma (H)=1$.
Let $K(H)=H^\perp_{N(X)}$ be the orthogonal complement to $H$ in $N(X)$.
Let $H^\ast=H/8$. Then any element $x\in N(X)$ can be written as
$$
x=aH^\ast+k^\ast
\tag{2.2.1}
$$
where $a\in \bz$ and $k^\ast\in K(H)^\ast$, because
$\bz H\oplus K(H)\subset N(X)\subset
N(X)^\ast\subset \bz H^\ast\oplus K(H)^\ast$. Since $\gamma(H)=1$,
the map $aH^\ast+[H] \mapsto k^\ast +K(H)$ gives an isomorphism of
the groups
$\bz/8\cong [H^\ast]/[H]\cong [u^\ast+K(H)]/K(H)$ where
$u^\ast+K(H)$ has order $8$ in $A_{K(H)}=K(H)^\ast/K(H)$.
It follows,
$$
N(X)=[\bz H, K(H), H^\ast+u^\ast].
\tag{2.2.2}
$$
The element $u^\ast$ is defined canonically mod $K(H)$. Since
$H^\ast+u^\ast$ belongs to the even lattice $N(X)$, it follows
$$
(H^\ast+u^\ast)^2={1\over 8}+{u^\ast}^2 \equiv 0 \mod 2.
\tag{2.2.3}
$$
Let $\overline{H^\ast}=H^\ast \mod [H]\in [H^\ast]/[H]\cong \bz/8$
and $\overline{k^\ast}=k^\ast\mod K(H)\in A_{K(H)}=K(H)^\ast/K(H)$.
We then have
$$
N(X)/[H,K(H)]=(\bz/8)(\overline{H^\ast} +\overline{u^\ast})\subset
(\bz/8)\overline{H^\ast}+K(H)^\ast/K(H).
\tag{2.2.4}
$$
Also $N(X)^\ast\subset \bz H^\ast+K(H)^\ast$ since
$H+K(H)\subset N(X)$,
and for $a\in \bz$, $k^\ast \in K(H)^\ast$ we have
$x=aH^\ast+k^\ast \in N(X)$ if and only if
$(aH^\ast+k^\ast)\cdot (H^\ast+u^\ast)={a\over 8}+k^\ast\cdot u^\ast\in \bz$.
It follows,
$$
N(X)^\ast=
\{aH^\ast+k^\ast\ |\ a\in \bz,\ k^\ast \in K(H)^\ast,\
a\equiv -8k^\ast\cdot u^\ast \mod 8 \}\subset \bz H^\ast+K(H)^\ast ,
\tag{2.2.5}
$$
and
$$
\split
N(X)^\ast/[H,K(H)]&=
\{(-8\overline{k^\ast}\cdot \overline{u^\ast})\, \overline{H^\ast}+
\overline{k^\ast}\}\ |\ \overline{k^\ast} \in A_{K(H)}\}
\subset \\
&\subset (\bz/8) \overline{H^\ast}+A_{K(H)}.
\endsplit
\tag{2.2.6}
$$
We introduce {\it the characteristic map of the polarization $H$}
$$
\kappa(H):K(H)^\ast \to A_{K(H)}/(\bz/8)(u^\ast+K(H))\to A_{N(X)}
\tag{2.2.7}
$$
where for $k^\ast \in K(H)^\ast$ we have
$$
\kappa(H)(k^\ast)=(-8k^\ast\cdot u^\ast) H^\ast+k^\ast + N(X)\in
A_{N(X)}.
\tag{2.2.8}
$$
It is epimorphic, its kernel is $(\bz/8)(u^\ast+K(H))$, and it gives
the canonical isomorphism
$$
\overline{\kappa(H)}:A_{K(H)}/(\bz/8)(u^\ast+K(H))
\cong A_{N(X)}.
\tag{2.2.9}
$$
For the corresponding discriminant forms we have
$$
\kappa(k^\ast)^2 \mod 2=(k^\ast)^2+8(k^\ast\cdot u^\ast )^2\mod 2.
\tag{2.2.10}
$$

Similar results we have for the polarization $h$ of $Y$.
We denote by $K(h)$ the orthogonal complement to $h$ in $N(Y)$.

\proclaim{Proposition 2.2.2} If $Y$ is isomorphic to $X$, then
the invariant $\gamma (H)$ of $H$ in $N(X)$ (see \thetag{1.1.1}) is
equal to $1$. If $\gamma (H)=1$ for $H\in N(X)$,
then $\gamma (h)=1$ for $h\in N(Y)$, and $\det K(h)=\det K(H)/4$.

Assume that $\gamma(H)=1$ for $H\in N(X)$.
Then the Mukai identification \thetag{2.1.5}
canonically identifies the transcendental
periods $(T(X), H^{2,0}(X))$ and
\newline
$(T(Y), H^{2,0}(Y))$. It follows that
the Picard lattices $N(Y)$ and $N(X)$ have the same genus.
In particular, $N(X)$ is isomorphic to $N(Y)$, if the genus of $N(Y)$
contains only one class. If the genus of $N(Y)$ contains only one class,
then $Y$ is isomorphic to $X$, if additionally the canonical homomorphism
$O(N(Y))\to O(q_{N(Y)})$ is epimorphic. Both these conditions are valid,
if either $\rho (Y)\ge 12$ or $N(Y)$ contains a copy of $U$
(see \thetag{1.1.2}).

Thus, $Y\cong X$, if $\gamma (H)=1$ for $H\in N(X)$ and
either $\rho(Y)\ge 12$ or $N(Y)$ contains a copy of $U$.
\endproclaim

Calculations below are valid for an arbitrary K3 surface $Y$ and
a primitive vector $h\in N(Y)$ with $h^2=2$ and $\gamma (h)=1$.
Let $K(h)$ be the orthogonal
complement to $h$ in $N(Y)$. Let $h^\ast=h/2$. Then any
element $x\in N(Y)$ can be written as
$$
x=ah^\ast+k^\ast
\tag{2.2.11}
$$
where $a\in \bz$ and $k^\ast\in K(h)^\ast$. Since $\gamma(h)=1$,
the map $ah^\ast+[h] \mapsto k^\ast +K(h)$ gives the isomorphism of
the groups
$\bz/2=[h^\ast]/[h]\cong [w^\ast+K(h)]/K(h)$ where
$w^\ast+K(h)$ has order $2$ in $K(h)^\ast/K(h)=A_{K(h)}$. It follows,
$$
N(Y)=[\bz h, K(h), h^\ast+w^\ast]
\tag{2.2.12}
$$
where $w^\ast+K(h)$ is an element of order $2$ in $A_{K(h)}$.
The element $w^\ast$ is defined canonically mod $K(h)$. The element
$h^\ast+w^\ast$ belongs to the even lattice $N(Y)$, it follows
$$
(h^\ast+w^\ast)^2={1\over 2}+{w^\ast}^2 \equiv 0 \mod 2.
\tag{2.2.13}
$$
Let $\overline{h^\ast}=h^\ast \mod [h]\in [h^\ast]/[h]\cong \bz/2$
and $\overline{k^\ast}=k^\ast\mod K(h)\in A_{K(h)}$.
We then have
$$
N(Y)/[h,K(h)]=(\bz/2)(\overline{h^\ast} +\overline{w^\ast})\subset
(\bz/2)\overline{h^\ast}+A_{K(h)}.
\tag{2.2.14}
$$
Since $N(Y)^\ast\subset \bz h^\ast+K(h)^\ast$,
for $a\in \bz$, $k^\ast \in K(h)^\ast$ we have
$x=ah^\ast+k^\ast \in N(Y)$ if and only if
$(ah^\ast+k^\ast)\cdot (h^\ast+w^\ast)={a\over 2}+k^\ast\cdot w^\ast\in \bz$.
It follows,
$$
N(Y)^\ast=
\{ah^\ast+k^\ast\ |\ a\in \bz,\ k^\ast \in K(h)^\ast,\
a\equiv -2k^\ast\cdot w^\ast \mod 2 \}\subset \bz h^\ast+K(h)^\ast \,,
\tag{2.2.15}
$$
and
$$
\split
N(Y)^\ast/[h,K(h)]&=
\{(-2\overline{k^\ast}\cdot \overline{w^\ast})\, \overline{h^\ast}+
\overline{k^\ast}\ |\ \overline{k^\ast} \in K(h)^\ast/K(h)\}
\subset \\
&\subset (\bz/2) \overline{h^\ast}+A_{K(h)}.
\endsplit
\tag{2.2.16}
$$
We introduce {\it the characteristic map of the polarization $h$}
$$
\kappa(h):K(h)^\ast \to A_{K(h)}/(\bz/2)(w^\ast+K(h))\to A_{N(Y)}
\tag{2.2.17}
$$
where for $k^\ast \in K(h)^\ast$ we have
$$
\kappa(h)(k^\ast)=(-2k^\ast\cdot w^\ast ) h^\ast+k^\ast + N(Y)\in
A_{N(Y)}.
\tag{2.2.18}
$$
It is epimorphic, its kernel is $(\bz/2)(w^\ast+K(h))$, and it gives
the canonical isomorphism
$$
\overline{\kappa(h)}:A_{K(h)}/(\bz/2)(w^\ast+K(h))=A_{N(Y)}.
\tag{2.2.19}
$$
For the corresponding discriminant forms we have
$$
\kappa(k^\ast)^2 \mod 2=(k^\ast)^2+2(k^\ast\cdot w^\ast )^2\mod 2.
\tag{2.2.20}
$$

Now we can formulate our main result:

\proclaim{Theorem 2.2.3} The surface $Y$ is isomorphic to $X$
if the following conditions (a), (b), (c) are valid:

(a) $\gamma(H)=1$ for $H\in N(X)$;

(b) there exists $\widetilde{h}\in N(X)$ with
$\widetilde{h}^2=2$, $\gamma(\widetilde{h})=1$ and such that
there exists an embedding
$f:K(H)\to K(\widetilde{h})$ of negative definite lattices such that
$K(\widetilde{h})=[f(K(H)), 4f(u^\ast)]$,
$w^\ast+K(\widetilde{h})=2f(u^\ast)+K(\widetilde{h})$;

(c) the dual to $f$ embedding
$f^\ast:K(\widetilde{h})^\ast \to K(H)^\ast$ commutes
(up to multiplication by $\pm 1$) with the characteristic maps
$\kappa(H)$ and $\kappa(\widetilde{h})$, i. e.
$$
\kappa(\widetilde{h})(k^\ast)=\pm \kappa(H)(f^\ast(k^\ast))
\tag{2.2.21}
$$
for any $k^\ast \in K(\widetilde{h})^\ast$.

The conditions (a), (b) and (c) are necessary if $\rk N(X)\le 19$ and
$X$ is a general K3 surface with the Picard lattice $N(X)$ in
the following sense: the automorphism group of the transcendental
periods $(T(X),H^{2,0}(X))$ is $\pm 1$. (Remind that,
by Proposition 2.2.1, $Y\cong X$ if $\rk N(X)=20$.)
\endproclaim

We can also formulate similar result using the surface $Y$.

\proclaim{Theorem 2.2.4} The surface $Y$ is isomorphic to $X$ if
the following conditions (a), (b) and (c) are valid:

(a) $\gamma (H)=1$ for $H\in N(X)$, then $\gamma (h)=1$ for $h\in N(Y)$;

(b) there exists a primitive
$\widetilde{H}\in N(Y)$ with
$\widetilde{H}^2=8$, $\gamma(\widetilde{H})=1$ and such that
there exists an embedding $f:K(\widetilde{H})\to K(h)$
of negative definite lattices such that
$K(h)=[f(K(\widetilde{H})), 4f(u^\ast)]$,
$w^\ast+K(h)=2f(u^\ast)+K(h)$;

(c) the dual to $f$ embedding
$f^\ast:K(h)^\ast \to K(\widetilde{H})^\ast$ commutes
(up to multiplication by $\pm 1$)
with the characteristic maps
$\kappa(\widetilde{H})$ and $\kappa(h)$, i. e.
$$
\kappa(h)(k^\ast)=\pm \kappa(\widetilde{H})(f^\ast(k^\ast))
\tag{2.2.22}
$$
for any $k^\ast \in K(h)^\ast$.

The conditions (a), (b) and (c) are necessary if $\rk N(Y)\le 19$ and
$Y$ is a general K3 surface with the Picard lattice $N(Y)$ in
the following sense:
the automorphism group of the transcendental periods $(T(Y),H^{2,0}(Y))$
is $\pm 1$. (Remind that, by Proposition 2.2.2,
$Y\cong X$ if $\gamma(H)=1$ and $\rk N(Y)=20$.)
\endproclaim

\subhead
2.3. Proofs
\endsubhead
Let us denote by $e_1$ the canonical generator of $H^0(X,\bz)$ and
by $e_2$ the canonical generator of $H^4(X,\bz)$. They generate
the sublattice $U$ in $H^\ast (X,\bz)$ with the Gram matrix $U$.
Consider Mukai vector $v=(2e_1+2e_2+H)$. We have
$$
N(Y)=v\,^\perp _{U\oplus N(X)}/\bz v.
\tag{2.3.1}
$$
Let us calculate $N(Y)$. Let $K(H)=(H)^\perp_{N(X)}$. Then we have
embedding of lattices of finite index
$$
\bz H\oplus K(H)\subset N(X)\subset N(X)^\ast\subset
\bz H^\ast \oplus K(H)^\ast
\tag{2.3.2}
$$
where $H^\ast =H/2$.
We have the orthogonal decomposition up to finite index
$$
U\oplus \bz H\oplus K(H)\subset U\oplus N(X)\subset
U\oplus \bz H^\ast \oplus K(H)^\ast .
\tag{2.3.3}
$$
Let $s=x_1e_1+x_2e_2+yH^\ast+z^\ast\in v^\perp_{U\oplus N(X)}$,
$z^\ast \in K(H)^\ast$.
Then $-2x_1-2x_2+y=0$ since $s\in v^\perp$ and hence $(s,v)=0$.
Thus, $y=2x_1+2x_2$ and
$$
s=x_1e_1+x_2e_2+2(x_1+x_2)H^\ast + z^\ast .
\tag{2.3.4}
$$
Here $s\in U\oplus N(X)$ if and only if $x_1, x_2 \in \bz$ and
$2(x_1+x_2)H^\ast+z^\ast\in N(X)$. This orthogonal complement
contains
$$
[\bz v, K(H), \bz h]
\tag{2.3.5}
$$
where $h=-e_1+e_2$, and this is a sublattice of finite index
in $(v^\perp)_{U\oplus N(X)}$. The generators $v$, generators of $K(H)$
and $h$ are free, and we can rewrite $s$ above using these generators
with rational coefficients as follows:
$$
s={-x_1+x_2\over 2}h+{x_1+x_2\over 4}v+z^\ast,
\tag{2.3.6}
$$
where $2(x_1+x_2)H^\ast+z^\ast \in N(X)$.
Equivalently,
$$
s=ah^\ast+b{v\over 4}+z^\ast,
\tag{2.3.7}
$$
where $a,b\in \bz$, $z^\ast \in K(H)^\ast$, $a\equiv b\mod 2$,
and $2bH^\ast+z^\ast\in N(X)$.

Thus, we get the following cases:

Assume that $\gamma(H)=8$. Then $2b\equiv 0\mod 8$, or
$b\equiv 0\mod 4$. Then $z^\ast \in K(H)$,
$a\equiv b\equiv 0\mod 2$, $s\in [h, K(H)]\mod \bz v$. It follows,
$$
N(X)=[H, K(H)]
\tag{2.3.8}
$$
and
$$
N(Y)=[h, K(h)=K(H)].
\tag{2.3.9}
$$
We have $\det N(Y)=\det N(X)/4$.

Assume that $\gamma(H)=4$. Then either $2b\equiv 0\mod 8$ or
$2b\equiv 4\mod 8$. Equivalently, $b\equiv 0,\ 2\mod 4$.
If $b\equiv 0\mod 4$, we get for $N(Y)$ the same elements as above.
If $b\equiv 2\mod 4$, we get an additional element $\mod \bz v$
which is equal to $z^\ast$ where $z^\ast\in K(H)^\ast$ is defined by the
condition that $4H^\ast +z^\ast \in N(X)$. Here $z^\ast+K(H)$ has
order two in $A_{K(H)}$. Thus, for this case
$$
N(X)=[H, K(H), {H\over 2}+z^\ast],
\tag{2.3.10}
$$
$$
N(Y)=[h, K(h)=[K(H),z^\ast]]
\tag{2.3.11}
$$
where $z^\ast+K(H)$ has order two in $K(H)^\ast/K(H)$.
We have $\det N(X)=2\det K(H)$ and $\det N(Y)=2\det K(H)/4=\det N(X)/4$.

Assume that $\gamma (H)=2$. Then $2b\equiv 0,\ 4,\ \pm 2\mod 8$. Or
$b\equiv 0,\ 2,\ \pm 1 \mod 4$. If $b\equiv 0, 2\mod 4$, we get the
same elements as for $\gamma (H)=4$. If $b\equiv \pm 1\mod 4$,
we get additional elements $\pm ({h\over 2}+z_1^\ast)$ where
$z_1^\ast+K(H)$ has order $4$ in $A_{K(H)}$. Finally we get
(changing notations) that
$$
N(X)=[H,K(H),{H\over 4}+z^\ast],
\tag{2.3.12}
$$
$$
N(Y)=[h, K(h)=[K(H),2z^\ast], {h\over 2}+z^\ast],
\tag{2.3.13}
$$
where $z^\ast+K(H)$ has order 4 in $K(H)^\ast/K(H)$.
We have $\det N(X)=\det K(H)/2$, $\det N(Y)=\det K(H)/8$,
and $\det N(Y)=\det N(X)/4$.

Assume that $\gamma (H)=1$. Then $2b\equiv 0,\ 4,\ \pm 2\mod 8$,
and we get the same lattice $N(Y)$ as above.
Thus,
$$
N(X)=[H,K(H),{H\over 8}+u^\ast],
\tag{2.3.14}
$$
$$
N(Y)=[h, K(h)=[K(H),4u^\ast],{h\over 2}+2u^\ast]=[h, K(h),{h\over 2}+w^\ast],
\tag{2.3.15}
$$
where $u^\ast+K(H)$ has order $8$ in $A_{K(H)}$, $w^\ast=2u^\ast$,
$K(h)=[K(H),2w^\ast=4u^\ast]$. Here we agreed notations with Sect. 2.2.
We have $\det N(X)=\det K(H)/8$ and $\det N(Y)=\det K(H)/8$.
Thus, $\det N(X)=\det N(Y)$ for this case.
We can formally put here $h={H\over 2}$ since
$h^2=\left({H\over 2}\right)^2=2$. Then
$$
N(X)\cap N(Y)=[H,K(H),{H\over 4}+2u^\ast].
\tag{2.3.16}
$$
We have
$$
[N(X):N(X)\cap N(Y)]=[N(Y):N(X)\cap N(Y)]=2.
\tag{2.3.17}
$$
From these calculations, we have:

\proclaim{Lemma 2.3.1} For Mukai identification \thetag{2.1.5},
the sublattice $T(X)\subset T(Y)$
has index $2$, if $\gamma(H)=2,\,4,\,8$, and
$T(X)=T(Y)$, if the $\gamma (H)=1$ for $H\in N(X)$.

Also $T(X)\subset T(Y)$ has index $2$ if either $\gamma(h)=2$
for $h\in N(Y)$
or $\det K(h)\not=\det K(H)/4$. If $\gamma (h)=1$ and
$\det K(h)=\det K(H)/4$, then $\gamma (H)=2$ or $1$, and
$T(X)\subset T(Y)$ has index $2$, if $\gamma (H)=2$, and
$T(X)=T(Y)$, if $\gamma (H)=1$ (we cannot get $T(X)=T(Y)$ using only the
Picard lattice $N(Y)$ of $Y$).
\endproclaim

\demo{Proof} Really, since $H\in N(X)$, $T(X)\perp N(X)$ and
$T(X)\cap \bz v=\{0\}$, the Mukai identification \thetag{2.1.5} gives
an embedding
$T(X)\subset T(Y)$. We then have $\det T(Y)=\det T(X)/[T(Y):T(X)]^2$.
Moreover, $|\det T(X)|=|\det N(X)|$ and $|\det T(Y)=|\det N(Y)|$
because the transcendental and the Picard lattice are orthogonal
complements to each other in a unimodular lattice $H^2(\ast,\bz)$.
By calculations above, we get the statement.

We remark that the first statement of Lemma 2.3.1 is a particular
case of the general statement by Mukai \cite{6} that $[T(Y):T(X)]=q$
where
$$
q=\min |v\cdot x|
\tag{2.3.18}
$$
for all $x\in H^0(X,\bz)\oplus N(X)\oplus H^4(X,\bz)$ such that
$v\cdot x\not=0$. For our Mukai vector $v=(2,\,H,\,2)$, it is easy to
see that $q=2$, if $\gamma(H)=2,\,4,\,8$, and $q=1$, if $\gamma (H)=1$.
\enddemo

Now we can prove Propositions 2.2.1 and 2.2.2. If $X\cong Y$, then
$T(X)\cong T(Y)$. Then $\det T(X)=\det T(Y)$, and $[T(Y):T(X)]=1$
for the Mukai identification. Then, $T(X)=T(Y)$ for the Mukai
identification \thetag{2.1.5}.
By Lemma 2.3.1, we then get first statements of Propositions 2.2.1 and 2.2.2.

Assume that $\gamma (H)=1$ as for Propositions 2.2.1 and 2.2.2.
Then $T(X)=T(Y)$ for the Mukai identification \thetag{2.1.5}.
By the discriminant forms technique (see \cite{8}),
then the discriminant
quadratic forms $q_{N(X)}=-q_{T(X)}$ and $q_{N(Y)}=-q_{T(Y)}$ are isomorphic.
Thus, lattices $N(X)$ and $N(Y)$ have the same signatures and
discriminant quadratic forms. It follows (see \cite{8}) that they have
the same genus: $N(X)\otimes \bz_p\cong N(Y)\otimes \bz_p$ for any prime
$p$ and the ring of $p$-adic integers $\bz_p$.
Additionally, assume that either the genus of $N(X)$ or the
genus of $N(Y)$ contains only one class. Then $N(X)$ and
$N(Y)$ are isomorphic.

If additionally the canonical homomorphism $O(N(X))\to O(q_{N(X)})$
(equivalently, $O(N(Y))\to O(q_{N(Y)})$) is epimorphic, then
the Mukai identification $T(X)=T(Y)$ can be extended to give an
isomorphism $\phi:H^2(X,\bz)\to H^2(Y,\bz)$ of cohomology lattices.
The Mukai identification is identical on $H^{2,0}(X)=H^{2,0}(Y)$.
Multiplying $\phi$ by $\pm 1$ and by elements of the reflection
group $W^{(-2)}(N(X))$, if necessary, we
can assume that $\phi(H^{2,0}(X))=H^{2,0}(Y)$ and $\phi$ maps the
K\"ahler cone of $X$ to the K\"ahler cone of $Y$. By global Torelli
Theorem for K3 surfaces proved by Piatetskii-Shapiro and
Shafarevich \cite{10}, $\phi$ is then defined by an isomorphism
of K3 surfaces $X$ and $Y$.

If $\rho (X)\ge 12$, by \cite{8, Theorem 1.14.4}, the primitive
embedding of $T(X)=T(Y)$ into the cohomology lattice $H^2(X,\bz)$ of
K3 surfaces is unique up to automorphisms of the lattice $H^2(X,\bz)$.
Like above, it then follows that $X$ is isomorphic to $Y$.

Let us prove Theorems 2.2.3 and 2.2.4.

Assume that $\gamma (H)=1$. The Mukai
identification then gives the canonical identification
$$
T(X)=T(Y).
\tag{2.3.19}
$$
Thus, it gives the canonical identifications
$$
\split
A_{N(X)}=N(X)^\ast/N(X)=&(U\oplus N(X))^\ast/(U\oplus N(X))=
T(X)^\ast/T(X)=A_{T(X)}\\
=A_{T(Y)}=T(Y)^\ast/T(Y)=&N(Y)^\ast /N(Y)=A_{N(Y)}.
\endsplit
\tag{2.3.20}
$$
Here $A_{N(X)}=N(X)^\ast/N(X)=(U\oplus N(X))^\ast/(U\oplus N(X))$ because
$U$ is unimodular,
$(U\oplus N(X))^\ast/(U\oplus N(X))=T(X)^\ast/T(X)=A_{T(X)}$ because
$U\oplus N(X)$ and $T(X)$ are orthogonal complements to each other in
the unimodular lattice $H^\ast(X,\bz)$. Here \linebreak
$A_{T(Y)}=T(Y)^\ast/T(Y)=N(Y)^\ast/ N(Y)=A_{N(Y)}$
because $T(Y)$ and $N(Y)$ are
orthogonal complements to each other in the unimodular lattice
$H^2(Y,\bz)$.
E. g. the identification
$(U\oplus N(X))^\ast/(U\oplus N(X))=T(X)^\ast/T(X)=A_{T(X)}$ is
given by the canonical correspondence
$$
x^\ast+(U\oplus N(X))\to t^\ast +T(X)
\tag{2.3.21}
$$
if $x^\ast \in (U\oplus N(X))^\ast$, $t^\ast \in T(X)^\ast$ and
$x^\ast+t^\ast \in H^\ast(X,\bz)$.

By \thetag{2.3.15}, we also have the canonical embedding of lattices
$$
K(H)\subset K(h)=[K(H),4u^\ast].
\tag{2.3.22}
$$

We have the key statement:

\proclaim{Lemma 2.3.2} The canonical embedding \thetag{2.3.22}
(it is given by \thetag{2.3.15}) $K(H)\subset K(h)$ of lattices,
and the canonical identification
$A_{N(X)}=A_{N(Y)}$ (given by \thetag{2.3.20}) agree with the
characteristic homomorphisms $\kappa(H):K(H)^\ast\to A_{N(X)}$
and $\kappa(h):K(h)^\ast \to A_{N(Y)}$, i.e.
$\kappa(h)(k^\ast)=\kappa (H)(k^\ast)$ for any $\kappa^\ast \in
K(h)^\ast\subset K(H)^\ast$ (this embedding is dual to \thetag{2.3.22}).
\endproclaim

\demo{Proof} From definitions of the identifications \thetag{2.3.20},
the identification \linebreak
$(U\oplus N(X))^\ast/(U\oplus N(X))=A_{N(Y)}$
is given by the canonical embeddings
$$
(U\oplus N(X))^\ast \supset (v^\perp )^\ast_0=(v^\perp/\bz v)^\ast
\tag{2.3.23}
$$
where
$(v^\perp)^\ast_0=\{s^\ast \in (U\oplus N(X))^\ast\ |\ s^\ast\cdot v=0\}$.

By \thetag{2.2.5},
$s^\ast=x_1e_1+x_2e_2+yH^\ast+k^\ast \in (U\oplus N(X))^\ast$
if and only if $x_1,x_2,y\in \bz$, $k^\ast \in K(H)^\ast$, and
$y\equiv -8k^\ast\cdot u^\ast\mod 8$. Here
$\kappa (k^\ast)=s^\ast+(U\oplus N(X))$.
We have $s^\ast \cdot v=-2x_1-2x_2+y$,
and $s^\ast \in (v^\perp)^\ast_0$, if additionally $y=2x_1+2x_2$.
Thus, we have
$2(x_1+x_2)\equiv -8k^\ast\cdot u^\ast \mod 8$. It follows,
$k^\ast \cdot (4u^\ast)\in \bz$, and
$k^\ast \in K(h)^\ast=[K(h), 4u^\ast]$. Moreover,
$x_1+x_2\equiv -2k^\ast \cdot w^\ast\mod 4$ and
$-x_1+x_2 \equiv x_1+x_2\equiv -2k^\ast\cdot  w^\ast\mod 2$
where $w^\ast=2u^\ast$.
Like in \thetag{2.3.7}, we then have that $s^\ast \in (v^\perp)^\ast_0$ if
and only if
$s^\ast =(-x_1+x_2)h^\ast+{x_1+x_2\over 4}v+k^\ast$ where $x_1$, $x_2$,
$k^\ast$ satisfy the conditions above. Finally, we have
$s^\ast \in (v^\perp)^\ast_0$, if and only if
$$
s^\ast =ah^\ast+b{v\over 4}+k^\ast
\tag{2.3.24}
$$
where $a, b\in \bz$, $k^\ast \in K(h)^\ast$,
$a\equiv -2k^\ast\cdot w^\ast\mod 2$
and $b\equiv -2k^\ast \cdot w^\ast\mod 4$. Here
$s^\ast$ gives $ah^\ast+k^\ast\in N(Y)^\ast$, and
$\kappa(h) (k^\ast)=ah^\ast+k^\ast+N(Y)=s^\ast+U\oplus N(X)=
\kappa (H)(k^\ast)$ under the identification
\thetag{2.3.20}. It proves the statement.
\enddemo

Proof of Theorem 2.2.3. We have the Mukai identification (it
is defined by \thetag{2.1.5}) of
the transcendental periods
$$
(T(X),\,H^{2,0}(X))=(T(Y),\,H^{2,0}(Y)).
\tag{2.3.25}
$$
For general $X$ with the Picard lattice $N(X)$, it is the unique
isomorphism of the transcendental periods up to multiplication by $\pm 1$.
If $X\cong Y$, this (up to $\pm 1$) isomorphism can be extended to
$\phi:H^2(X, \bz)\cong H^2(Y, \bz)$. The restriction of $\phi$ on $N(X)$
gives then isomorphism $\phi_1:N(X)\cong N(Y)$ which is $\pm 1$ on
$A_{N(X)}=A_{N(Y)}$ under the identification \thetag{2.3.20}.
The element $\widetilde{h}=(\phi_1)^{-1}(h)$ and $f=\phi^{-1}$
satisfy Theorem 2.2.3 by Lemma 2.3.2.

The other way round, under conditions of Theorem 2.2.3, by Lemma 2.3.2,
one can construct an isomorphism $\phi_1:N(X)\cong N(Y)$ which is
$\pm 1$ on $A_{N(X)}=A_{N(Y)}$. It can be extended to
be $\pm 1$ on the transcendental periods under the Mukai identification
\thetag{2.3.25}. Then it is defined by the isomorphism
$\phi:H^2(X,\,\bz)\to H^2(Y,\,\bz)$. Multiplying $\phi$ by $\pm 1$ and
by reflections from $W^{(-2)}(N(X))$, if necessary (the group
$W^{(-2)}(N(X))$ acts identically on the discriminant group
$N(X)^\ast/N(X)$), we can assume that $\phi$ maps the
K\"ahler cone of $X$ to the K\"ahler cone of $Y$.
By global Torelli Theorem \cite{10}, it is then defined by an isomorphism
of $X$ and $Y$.

To prove Theorem 2.2.4, one should argue similarly.

\head
3. The case of Picard number 2
\endhead

\subhead
3.1. General results
\endsubhead
Here we want to apply results of Sect. 2 to $X$ and $Y$ with
Picard number 2.

We start with some preliminary considerations on
K3 surfaces with Picard number 2 and primitive polarizations of
degree 8 and 2.

Assume that $\rk N(X)=2$. Let $H\in N(X)$ be a primitive
polarization with $H^2=8$. Assume that $\gamma (H)=1$ for $H\in N(X)$
(we have this condition, if $Y\cong X$). Let
$K(H)=H^\perp_{N(X)}=\bz \delta$ and $\delta^2=-t$ where $t>0$ is even.
It then follows that
$N(X)=[\bz H,\, \bz \delta,\, H^\ast+\mu {\delta\over 8}]$
where $H^\ast=H/8$ and $\mu=\pm 1,\,\pm 3$. Since
$(H^\ast+\mu {\delta\over 8})^2={1/8-\mu^2t/64}\equiv 0\mod 2$,
it follows  $t=8d$ where $d\in \bn$ and $1-\mu^2d\equiv 0\mod 16$.
Then $d$ is odd and $\mu^2d\equiv 1\mod 16$. It follows that
$d\equiv 1\mod 8$ and
$\mu^2\equiv d\mod 16$. Changing $\delta$ by $-\delta$ if
necessary, we can always assume that $\mu=1$ or $\mu=3$ where
$\mu =1$ if $d\equiv 1\mod 16$, and $\mu=3$ if $d\equiv 9\mod 16$.
These simple calculations show that the lattice $N(X)$ and
$H$ are defined uniquely, up to isomorphisms, by $d$ where
$-d=\det N(X)$. One can even replace $H$ by
any primitive element of $N(X)$ with square $8$. Thus, we have

\proclaim{Proposition 3.1.1} Let $X$ be a K3 surface with the Picard
number $\rho =2$. Assume that $X$ has a primitive
polarization $H$ of degree $H^2=8$, and $\gamma (H)=1$ for $H\in N(X)$.
Then the lattice $N(X)$ is defined by its
determinant $\det N(X)=-d$ where $d\in \bn$ and $d\equiv 1\mod 8$.
There exists a unique choice of a primitive orthogonal
vector $\delta \in K(H)=H^\perp_{N(X)}$ such that
$$
N(X)=[H,\delta,(H+\mu\delta)/8]
\tag{3.1.1}
$$
where $\delta^2=-8d$,
$$
\mu=
\cases
1,\ &\text{if $d\equiv 1\mod 16$},\\
3,\ &\text{if $d\equiv 9\mod 16$}.
\endcases
\tag{3.1.2}
$$
We have
$$
N(X)=\{z=(xH+y\delta)/8\ |\ x,y\in \bz\ \text{and}\ \mu x\equiv y\mod 8\}.
\tag{3.1.3}
$$

For any primitive element $H^\prime \in N(X)$ with
$(H^\prime)^2 =H^2=8$, there exists a unique automorphism
$\phi\in O(N(X))$ such that $\phi(H)=H^\prime$.

We denote the above (unique up to isomorphisms)
hyperbolic lattice $N(X)$ by $N^{8}_d$
where $d\in \bn$ and $d\equiv 1\mod 8$.
\endproclaim

\demo{Proof} Let $H^\prime$ be a primitive element of $N(X)\cong N^{8}_d$
with square $8$. It is easy to see that, if
$\gamma (H^\prime)\not=1$, then $\det N(X)$ is even which is impossible.
The rest calculations are elementary.
\enddemo

Now let $Y$ be a K3 surface with $\rk N(Y)=2$.
Let $h\in N(Y)$ be a $nef$ element of degree $h^2=2$.
Assume that $\gamma (h)=1$ (this condition is necessary to have
$Y\cong X$). Let $K(h)=h^\perp_{N(Y)}=\bz\alpha$ and
$\alpha^2=-t$ where $t>0$ is even.
It then follows that $N(Y)=[\bz h,\,\bz \alpha,\,
h^\ast +{\alpha \over 2}]$ where $h^\ast =h/2$. We have
$(h^\ast+{\alpha\over 2})^2=1/2-t/4\equiv 0\mod 2$, and hence
$t=2d$ where $d\in \bn$ and $d\equiv 1\mod 4$. Like for $X$ above,
we get

\proclaim{Proposition 3.1.2} Let $Y$ be a K3 surface with the Picard
number $\rho =2$. Assume that $Y$ has a $nef$ element
$h\in N(Y)$ of degree $h^2=2$, and $\gamma (h)=1$ for
$h\in N(Y)$. Then the lattice $N(Y)$ is defined by its
determinant $\det N(Y)=-d$ where $d\in \bn$ and $d\equiv 1\mod 4$.
For a primitive orthogonal
vector $\alpha \in K(h)=h^\perp_{N(Y)}$
(it is unique up to changing by $-\alpha$), we have
$$
N(Y)=[h,\,\alpha,\,(h+\alpha)/2]
\tag{3.1.4}
$$
where $\alpha^2=-2d$.
We have
$$
N(Y)=\{z=(xh+y\alpha)/2\ |\ x,y\in \bz\ \text{and}\ x\equiv y\mod 2\}.
\tag{3.1.5}
$$

For any element $h^\prime \in N(Y)$ with
$(h^\prime)^2 =h^2=2$, there exists an automorphism
$\phi\in O(N(Y))$ such that $\phi(h)=h^\prime$.
It is unique up to changing to $\phi s_\alpha$ where
$s_\alpha(h)=h$, $s_\alpha(\alpha)=-\alpha$.

We denote the above (unique up to isomorphisms)
hyperbolic lattice $N(Y)$ of the determinant $-d$ by $N^{2}_d$
where $d\in \bn$ and $d\equiv 1\mod 4$.
\endproclaim

\demo{Proof} It is trivial.
\enddemo

From Propositions 3.1.1 and 3.1.2, we then get

\proclaim{Proposition 3.1.3} Under conditions and notations of
Propositions 3.1.1 and 3.1.2, all elements
$h^\prime=(xH+y\delta)/8\in N(X)$ with
$(h^\prime)^2=2$ are in one to one
correspondence with solutions $(x,y)$ of the equation
$$
x^2-dy^2=16
\tag{3.1.6}
$$
with odd $x$, $y$ and $\mu x\equiv y\mod 8$. Changing an odd
solution $(x,y)$ of the equation to $(x,-y)$, if necessary,
one can always satisfy the last congruence.

The Picard lattices of $X$ and $Y$ are isomorphic, $N(X)\cong N(Y)$,
if and only if $\det N(X)=\det N(Y)=-d$ (it follows that
$d\equiv 1\mod 8$), and there exists an element
$h^\prime\in N(X)$ with $(h^\prime)^2=2$ as above.
Equivalently, the equation $x^2-dy^2=16$ has a solution
with odd $x$ and $y$.
\endproclaim

\demo{Proof}
Assume, $h^\prime=(xH+y\delta)/8\in N(X)$ and
$(h^\prime)^2=2$. This is equivalent to $x,y\in \bz$,
$\mu x\equiv y\mod 8$ and $2=(x^2-dy^2)/8$. Thus, the elements $h^\prime$
are in one to one correspondence with integral solutions $(x,y)$ of
the equation $x^2-dy^2=16$ which satisfy the condition
$x\mu \equiv y\mod 8$.

Let $(x,y)$ be an integral solution of the equation. Clearly, then
$x\equiv y\mod 2$. Let $(x,y)$ be even. Then $(x,y)=2(x_1,y_1)$ where
$(x_1,y_1)$ is an integral solution of the equation $x_1^2-dy_1^2=4$.
Again $x_1\equiv y_1\mod 2$. If $x_1\equiv y_1\equiv 1\mod 2$,
we get a contradiction with $d\equiv 1\mod 8$. If
$x_1\equiv y_1\equiv 0\mod 2$, we get that $(x,y)=4(x_2,y_2)$ where
$(x_2,y_2)$ are integral solutions of the equation $x_2^2-dy_2^2=1$.
It follows that $x_2\mod 2$ and $y_2\mod 2$ are different.
Then the congruence $\mu x\equiv y\mod 8$ is not satisfied, and the
solution $(x,y)$ does not give an element of $N(X)$.

Assume that $(x,y)$ is a solution of $x^2-dy^2=16$ with odd $x$, $y$.
We have $d\equiv \mu^2\mod 16$ and $(x-y\mu )(x+y\mu )\equiv 0\mod 16$.
If $x-y\mu \equiv 0\mod 4$ and $x+y\mu\equiv 0\mod 4$, then
$2x\equiv 0\mod 4$ which is impossible for odd $x$. Thus, we
have that only one of congruences: $x-y\mu\equiv 0\mod 8$ or
$x+y\mu\equiv 0\mod 8$ is valid.
It follows that exactly one of solutions $(x,y)$ or
$(x,-y)$ gives an element of $N(X)$. It proves first statement.

If $N(Y)\cong N(X)$, the lattices have the same
determinant. By Proposition 3.1.2, $N(X)\cong N(Y)$, if and only if $N(X)$
has an element $h^\prime$ with square $2$. This finishes the proof.
\enddemo

We have a similar statement in terms of $Y$.

\proclaim{Proposition 3.1.4} Under conditions and notations of
Propositions 3.1.1 and 3.1.2,
all primitive elements
$H^\prime=(xh+y\alpha)/2\in N(Y)$ with square $(H^\prime)^2=8$
are in one to one correspondence with solutions $(x,y)$ of the equation
$$
x^2-dy^2=16
$$
with odd $x$, $y$.

The Picard lattices of $X$ and $Y$ are isomorphic, $N(X)\cong N(Y)$,
if and only if $\det N(X)=\det N(Y)=-d$ (it follows that
$d\equiv 1\mod 8$), and there exists an element $H^\prime\in N(Y)$ with
$(H^\prime)^2=8$ as above.
Equivalently, the equation $x^2-dy^2=16$ has a solution with
odd $x$ and $y$. It follows that $d\equiv 1\mod 8$.
\endproclaim

\demo{Proof}
Let $H^\prime=(xh+y\alpha)/2$ where $x,y\in \bz$ and $x\equiv y\mod 2$.
Then $2=(x^2-dy^2)/8$, hence $x^2-dy^2=16$.
The element $H^\prime$ is not primitive if and only if
$x/2,y/2\in \bz$ and $x/2\equiv y/2\mod 2$. This is equivalent that
$x\equiv y\equiv 0\mod 2$ and $x\equiv y\mod 4$.
Thus, the elements $H^\prime$
are in one to one correspondence with integral solutions $(x,y)$ of
the equation $x^2-dy^2=16$ which satisfy the condition that
either they are both odd or both even, but
$x\mod 4$ and $y\mod 4$ are different.

Assume that $(x,y)$ is an integral solution of the equation
$x^2-dy^2=16$. Clearly $x\equiv y\mod 2$ since $d\equiv 1\mod 8$.
If $x\equiv y\equiv 1\mod 2$, then $(x,y)$ gives a primitive element 
$H^\prime$ by previous considerations.
Let $(x,y)$ be even. Then $(x,y)=2(x_1,y_1)$ where
$(x_1,y_1)$ is an integral solution of the equation $x_1^2-dy_1^2=4$.
Again $x_1\equiv y_1\mod 2$. It then follows that $x\equiv y\mod 4$,
and the solution $(x,y)$ does not give a primitive element
$H^\prime$ of $N(Y)$.

If $N(Y)\cong N(X)$, the lattices have the same
determinant. By Proposition 3.1.1, $N(X)\cong N(Y)$,
if and only if the lattice
$N(Y)$ has a primitive element $H^\prime$ with
$(H^\prime)^2=8$. It finishes the proof.
\enddemo

Now we can apply Theorems 2.2.3 and 2.2.4 to find out when
$Y\cong X$.

\proclaim{Theorem 3.1.5}
Let $X$ be a K3 surface and $\rho (X)=2$.
Assume that $X$ is an intersection of three quadrics (more generally,
$X$ has a primitive polarization $H$ of degree 8). Let $Y$ be a K3 surface
which is the double covering of the net $\bp^2$ of the quadrics defining
$X$, ramified along the curve of degenerate quadrics and $h$ be the
preimage of a line in $\bp^2$ (more generally, $Y$ is the moduli space
of sheaves on $X$ with the Mukai vector $v=(2,H,2)$ and the canonical
$nef$ element $h=(-1,0,1)\mod \bz v$).

If $Y\cong X$, then
$$
\gamma(H)=1,\ \det N(X)=-d\ \text{where\ } d\equiv 1\mod 8.
\tag{3.1.7}
$$
In the case  \thetag{3.1.7}, with notations of Propositions 3.1.1 and
3.1.3 for $X$ and $H$,
we have that all elements $\widetilde{h}=(xH+y\delta)/8\in N(X)$ with square
$\widetilde{h}^2=2$ satisfying Theorem 2.2.3 are in one to one correspondence
with integral solutions $(x,y)$ of the equation
$$
x^2-dy^2=16
\tag{3.1.8}
$$
with odd $x,\,y$, and $x\equiv \pm 4\mod d$,
$\mu x\equiv y\mod 8$ (one can always satisfy the last
congruence changing $y$ to $-y$ if necessary).

In particular (by Theorem 2.2.3),
for a general K3 surface $X$ with $\rho(X)=2$ and
$\det N(X)=-d$, $d\equiv 1\mod 8$, we have:
$Y\cong X$ if and only if the equation $x^2-dy^2=16$ has an integral
solution $(x,y)$ with odd $(x,y)$ and $x\equiv \pm 4\mod d$.
Moreover, a nef element $h=(xH+y\delta)/8$ of $X$ with $h^2=2$ defines
the structure of a double plane on $X$ which is isomorphic to the
double plane $Y$ if and only if $x\equiv \pm 4\mod d$.
\endproclaim

\demo{Proof} If $Y\cong X$, then $\gamma (H)=1$ by Proposition 2.2.1.
By Proposition 3.1.1, $N(X)\cong N^8_d$ where $d\in \bn$,
$d\equiv 1\mod 8$ and $\det N(X)=-d$.

Assume that $Y\cong X$ for a general $X$ with the Picard lattice $N^8_d$.
Let $\widetilde{h}\in N(X)$ satisfies conditions of Theorem 2.2.3.

By Proposition 3.1.3, all primitive
$$
\widetilde{h}=(xH+y\delta)/8
\tag{3.1.9}
$$
with $(\widetilde{h})^2=2$ are in one to one
correspondence with odd $(x,y)$ which satisfy the equation
$x^2-dy^2=16$ and $y\equiv \mu x\mod 8$,
and any integral solution of the equation $x^2-dy^2=16$ with odd
$x,\,y$ gives such $\widetilde{h}$
after replacing $y$ to $-y$ if necessary,
which does not matter for the statement of Theorem.

Let $k=aH+b\delta \in \widetilde{h}^\perp =\bz \alpha$. Then
$(k,h)=ax-byd=0$ and $(a,b)=\lambda(yd,x)$. Hence, we have
$(\lambda(ydH+x\delta))^2=\lambda^2(8y^2d^2-8dx^2)=8\lambda^2d(y^2d-x^2)=
-2^7\lambda^2d$. Since $\alpha^2=-2d$, we get $\lambda=2^{-3}$ and
$\alpha=(ydH+x\delta)/8$. There exists a unique (up to $\pm 1$)
embedding $f:K(H)=\bz\delta\to K(h)=\bz\alpha$
of one-dimensional lattices. It is given by $f(\delta)=2\alpha$
up to $\pm 1$. Thus, its dual is defined by
$f^\ast(\alpha^\ast)=2\delta^\ast$ where $\alpha^\ast=\alpha/2d$ and
$\delta^\ast=\delta/8d$. To satisfy
conditions of Theorem 2.2.3, we should have
$$
\kappa(\widetilde{h})(\alpha^\ast)=\pm 2\kappa(H)(\delta^\ast).
\tag{3.1.10}
$$
We have $u^\ast=\mu d\delta^\ast$, $w^\ast=2f(u^\ast)=\mu d\alpha^\ast=
\mu \alpha/2$, and
$$
\kappa(\widetilde{h})(\alpha^\ast)=(-2\alpha^\ast\cdot w^\ast)
\widetilde{h}^\ast+
\alpha^\ast+N(X)
\tag{3.1.11}
$$
by \thetag{2.2.18}. Here $\widetilde{h}^\ast=\widetilde{h}/2$.
We then have $\alpha^\ast\cdot w^\ast=\mu (\alpha/2d)\cdot (\alpha/2)=
-\mu/2$, and
$\kappa(\widetilde{h})(\alpha^\ast)=\mu \widetilde{h}^\ast+\alpha^\ast+N(X)=
\mu (x(H/16)+y(\delta/16))+y(H/16)+x\delta/16d$. It follows,
$$
\kappa(\widetilde{h})(\alpha^\ast)={\mu x+y\over 16}H+{x+\mu yd\over 16d}
\delta+N(X)=
{\mu x+y\over 2}H^\ast+{x+\mu yd\over 2}\delta^\ast+N(X) 
\tag{3.1.12}
$$
where $H^\ast=H/8$. We have
$u^\ast=\mu d\delta^\ast=\mu \delta/8$.
By \thetag{2.2.8},
$\kappa(H)(\delta^\ast)=
(-8\delta^\ast\cdot u^\ast)H^\ast+\delta^\ast+N(X)$.
We have $\delta^\ast\cdot u^\ast=-\mu/8$. It follows,
$$
\kappa (H)(\delta^\ast)=\mu H^\ast+\delta^\ast.
\tag{3.1.13}
$$
By \thetag{3.1.12} and \thetag{3.1.13}, we then get that
$\kappa (H)(2\delta^\ast)=\pm \kappa (\widetilde{h})(\alpha^\ast)$ is
equivalent to $(x+\mu yd)/2\equiv \pm 2\mod d$ or
$x+\mu yd\equiv \pm 4\mod d$ since
the group $N(X)^\ast/N(X)$ is cyclic of order $d$ and it is generated by
$\mu H^\ast+\delta^\ast+N(X)$. Thus, finally we get
$x\equiv \pm 4\mod d$.

The condition to have an odd solution $(x,y)$ of \thetag{3.1.18} with
$x\equiv \pm 4\mod d$ does not depend on the choice of $X$ and
the polarization $H$. If this condition is
satisfied, it is valid for any $X$ with the Picard lattice $N(X)$ and
any its primitive polarization $H$ of degree 8.
By Theorem 2.2.3, $Y\cong X$ for all of them.
This finishes the proof.
\enddemo

We have a similar statement in terms of $Y$.

\proclaim{Theorem 3.1.6}
Let $X$ be a K3 surface and $\rho (X)=2$.
Assume that $X$ is an intersection of three quadrics (more generally,
$X$ has a primitive polarization $H$ of degree 8).
Let $Y$ be a K3 surface which is the double covering
of the net $\bp^2$ of the quadrics
ramified along the curve of degenerate quadrics and
$h\in N(Y)$ be the preimage of
line in $\bp^2$ (more generally, $Y$ is the moduli space of sheaves
on $X$ with the Mukai vector $v=(2,H,2)$ and the canonical
$nef$ element $h=(-1,0,1)\mod \bz v$).

If $Y\cong X$, then
$$
\gamma (H)=1,\ \gamma(h)=1,\ \det N(Y)=-d\ \text{where\ } d\equiv 1\mod 8.
\tag{3.1.14}
$$
In the case \thetag{3.1.14} with notations of
Propositions 3.1.2 and 3.1.4 for $Y$ and $h$, then we have that
all elements $\widetilde{H}=(xh+y\alpha)/2\in N(Y)$ with square
$\widetilde{H}^2=8$ satisfying Theorem 2.2.4 are in one to one correspondence
with integral solutions $(x,y)$ of the equation
$$
x^2-dy^2=16
\tag{3.1.15}
$$
with odd $x,\,y$, and $x\equiv \pm 4\mod d$.

In particular (by Theorem 2.2.4),
for a general K3 surface $Y$ with $\rho(Y)=2$ and
$\det N(Y)=-d$, $d\equiv 1\mod 8$, we have:
$Y\cong X$ if and only if the equation $x^2-dy^2=16$ has an integral
solution $(x,y)$ with odd $(x,y)$ and $x\equiv \pm 4\mod d$. Moreover,
a primitive polarization $H=(xh+y\alpha)/2$ of $Y$ with
$H^2=8$ defines a structure of intersection of three quadrics on $Y$
which is isomorphic to the structure of intersection of three quadrics
of $X$ if and only if $x\equiv \pm 4\mod d$.
\endproclaim

\demo{Proof} If $Y\cong X$, then $\gamma (H)=1$ and
$\gamma (h)=1$, by Proposition 2.2.2.
By Propositions 3.1.2 and 3.1.4, $N(Y)\cong N^2_d$ where $d\equiv 1\mod 8$
and $\det N(Y)=-d$.

Assume that $Y\cong X$ for a general $Y$ with the Picard lattice $N^2_d$.
Let $\widetilde{H}\in N(Y)$ satisfies conditions of Theorem 2.2.4.

By Proposition 3.1.4, all primitive
$$
\widetilde{H}=(xh+y\alpha)/2
\tag{3.1.16}
$$
with $\widetilde{H}^2=8$
are in one to one correspondence with odd $(x,y)$ which
satisfy the equation $x^2-dy^2=16$.

Let $k=ah+b\alpha \in \widetilde{H}^\perp_{N(Y)} =\bz \delta$. Then
$ax-byd=0$ and $(a,b)=\lambda(yd,x)$. We have
$(\lambda(ydh+x\alpha))^2=\lambda^2(2y^2d^2-2dx^2)=
2\lambda^2d(y^2d-x^2)=
-2^5\lambda^2d$. Since $\delta^2=-8d$, we get $\lambda=2^{-1}$, and
$$
\delta=(ydh+x\alpha)/2.
\tag{3.1.17}
$$
There exists a unique (up to $\pm 1$) embedding
$f:K(\widetilde{H})=\bz\delta\to K(h)=\bz\alpha$
of one-dimensional lattices. It is given by $f(\delta)=2\alpha$
up to $\pm 1$. Thus, its dual is defined by
$f^\ast(\alpha^\ast)=2\delta^\ast$ where $\alpha^\ast=\alpha/2d$ and
$\delta^\ast=\delta/8d$. To satisfy
conditions of Theorem 2.2.4, we should have
$$
\kappa(h)(\alpha^\ast)=\pm 2\kappa(\widetilde{H})(\delta^\ast).
\tag{3.1.18}
$$
Like for the proof above, we have
$$
\kappa(h)(\alpha^\ast)=\mu h^\ast+\alpha^\ast+N(Y)
\tag{3.1.19}
$$
where $h^\ast=h/2$, and
$$
\kappa (\widetilde{H})(\delta^\ast)=\mu \widetilde{H}^\ast+\delta^\ast
\tag{3.1.20}
$$
where $\widetilde{H}^\ast=\widetilde{H}/8$.
By \thetag{3.1.16} and \thetag{3.1.17}, we then have
$$
\kappa (\widetilde{H})(\delta^\ast)=\mu (xh+y\alpha)/16+(ydh+x\alpha)/16d+N(Y).
\tag{3.1.21}
$$
It follows,
$$
\kappa (\widetilde{H})(\delta^\ast)={\mu x+y\over 8}h^\ast+{x+\mu dy\over 8}
\alpha^\ast+N(Y).
\tag{3.1.22}
$$
Thus,
$\kappa (\widetilde{H})(2\delta^\ast)=\pm \kappa (h)(\alpha^\ast)$
is equivalent to $(x+\mu yd)/4\equiv \pm 1\mod d$ or
$x+\mu yd\equiv \pm 4\mod d$ since
the group $N(Y)^\ast/N(Y)$ is cyclic of the order $d$, and it is generated by
$\mu h^\ast+\alpha^\ast+N(Y)$. Thus, finally we get
$x\equiv \pm 4\mod d$.

The condition to have an odd solution $(x,y)$ of \thetag{3.1.15} with
$x\equiv \pm 4\mod d$ does not depend on the choice of $Y$ and
the polarization $h$. If this condition is
satisfied, it is valid for any $Y$ with the Picard lattice $N(Y)$ and
any its polarization $h$ of degree 2.
By Theorem 2.2.4, $Y\cong X$ for all of them.
This finishes the proof.
\enddemo

By Theorems 3.1.5 and 3.1.6, general $X$ and $Y$ with $\rho=2$ and
$X\cong Y$ are labeled by $d\in \bn$ such that the conditions
$$
x^2-dy^2=16,\ x\equiv y\equiv 1 \mod 2,\ x\equiv \pm 4\mod d
\tag{3.1.23}
$$
are satisfied for some integers $x$, $y$.
Since $x$ and $y$ are odd, it follows that
$d\equiv 1\mod 8$, which we had known. Since $x\equiv \pm 4\mod d$
and $x$ is odd, we can write down all these $x$ as
$$
x=\pm 4+kd,\ k\in \bz,\ k\equiv 1\mod 2.
\tag{3.1.24}
$$
The equation \thetag{3.1.23} gives then
$16\pm 8kd+k^2d^2-dy^2=16$, i. e.  $\pm 8k+k^2d-y^2=0$.
It follows
$$
d={y^2+8k\over k^2}
\tag{3.1.25}
$$
where $k$ is any odd integer. Let $p$ be an odd prime.
Assume that $p^{2t+1}|k$, but
$p^{2t+2}\nmid k$ for some $t\ge 0$.
Since $k|(y^2+8k)$, it follows $k|y^2$ and $p^{t+1}|y$.
Then $p^{2t+2}|y^2$.
Since $p^{4t+2}|k^2$, it then follows that $p^{2t+2}|y^2+8k$ and
$p^{2t+2}|y^2$. Thus, $p^{2t+2}|k$. We get a contradiction.
It shows that an odd prime $p$ may divide $k$ only in an
even maximal power. Thus $k=\mp b^2$ for some $b\in \bn$.
Then $b|y$ (since $k|y^2$), and $y=ab$ for some integer $a$.
Thus, finally we get that
$$
d={a^2\mp 8\over b^2}
\tag{3.1.26}
$$
for some odd integers $a$ and $b$. Then $(a,b)$ is an
odd solution of one of equations
$$
a^2-db^2=8,
\tag{3.1.27}
$$
or
$$
a^2-db^2=-8.
\tag{3.1.28}
$$
If one of these equations has a solution with odd $a$ and $b$,
then $d\equiv 1\mod 8$. Further we assume that: $d\equiv 1\mod 8$.
It is easy to see that then any solution $(a,b)$ of
the equations \thetag{3.1.27} and
\thetag{3.1.28} has odd $a$
and $b$. Our considerations show that solutions $(x,y)$ of
\thetag{3.1.23} are in one to one correspondence with solutions
$(a,b)$ of
the equations \thetag{3.1.27} and \thetag{3.1.28}.
More exactly, we get solutions
$$
(x,y)=\pm \left(a^2-4,\ ab\right),\ \text{if\ }a^2-db^2=8,
\tag{3.1.29}
$$
and
$$
(x,y)=\pm \left(a^2+4,\ ab\right),\ \text{if\ }a^2-db^2=-8
\tag{3.1.30}
$$
of \thetag{3.1.23}, and any solution of \thetag{3.1.23} can be written in
this form. We call solutions \thetag{3.1.29} and \thetag{3.1.30}
of \thetag{3.1.23} as {\it associated with solutions
$(a,b)$ of the equations \thetag{3.1.27} and \thetag{3.1.28} respectively.}

Thus, we get the final result

\proclaim{Theorem 3.1.7}
Let $X$ be a K3 surface and $\rho (X)=2$.
Assume that $X$ is an intersection of three quadrics (more generally,
$X$ has a primitive polarization $H$ of degree 8). Let $Y$ be a K3 surface
which is the double covering of the net $\bp^2$ of the quadrics defining
$X$ ramified along the curve of degenerate quadrics and $h$ be preimage of a
line in $\bp^2$ (more generally, $Y$ is the moduli space of sheaves
on $X$ with the Mukai vector $v=(2,H,2)$ and the canonical
$nef$ element  $h=(-1,0,1)\mod \bz v$).

Then $Y\cong X$ for a general $X$ with $\rho (X)=2$, if and only if
$\det N(X)=-d$ where $d\equiv 1\mod 8$, and one of equations
$$
a^2-db^2=8,
\tag{3.1.31}
$$
or
$$
a^2-db^2=-8
\tag{3.1.32}
$$
has an integral solutions. All solutions of these equations have
odd $a$ and $b$, and the set of possible $d$ is union of the
two infinite sets
$$
\Da_+ =\{{a^2-8\over b^2}\in \bn\ |\ a, b\in \bn\
\text{are odd}\}
\tag{3.1.33}
$$
for the equation \thetag{3.1.31}, and
$$
\Da_- =
\{{a^2+8\over b^2}\in \bn\ |\ a, b\in \bn\ \text{are odd}\}
\tag{3.1.34}
$$
for the equation \thetag{3.1.32}. Here $\Da_+$ is infinite because it
contains the infinite subset $\{a^2-8 |\ a\in \bn\
\text{is odd}\}$, and $\Da_-$ is infinite because it contains
the infinite subset $\{a^2+8 |\ a\in \bn\
\text{is odd}\}$.

Solutions of \thetag{3.1.31} and \thetag{3.1.32} give all solutions of
\thetag{3.1.23}
as associated solutions \thetag{3.1.29} and \thetag{3.1.30}, and
all primitive elements
$\widetilde{h}\in N(X)$ with $\widetilde{h}^2=2$ of Theorem 3.1.5, and
$\widetilde{H}\in N(Y)$ with $\widetilde{H}^2=8$ of Theorem 3.1.6.
\endproclaim

We also mention the following: {\it The sets
$\Da_+-\Da_+\cap \Da_-$ and $D_--\Da_+\cap \Da_-$ are infinite.}
The infinite sequence $\{(1+14k)^2-8\ |\ k\in\bn\}$ is in
$\Da_+ - \Da_+\cap \Da_-$. Its elements $d$ are not in
$\Da_-$ because $7\vert d$ and $-2$ is not a square
$\mod 7$. Similarly, the infinite sequence
$\{(1+6k)^2 + 8\ |\ \ k\in\bn\}$ is in $D_--\Da_+\cap \Da_-$.
We don't know: {\it If the set $\Da_+\cap \Da_-$ is infinite?}

\smallpagebreak

It is interesting to interpret solutions of
\thetag{3.1.31} and \thetag{3.1.32} as appropriate elements of
$N(X)$ and $N(Y)$.

We have

\proclaim{Theorem 3.1.8} Under conditions and notations of
Proposition 3.1.1, the elements
$$
h_1=(2aH+2b\delta)/8\in N(X),
\tag{3.1.35}
$$
where $(a,b)$ is any integral solution of $a^2-db^2=\pm 8$ satisfying
the congruence $\mu a\equiv b\mod 4$ (one can always satisfy the
congruence changing $b$ to $-b$ if necessary) are all elements
$h_1\in N(X)$ with
$$
(h_1)^2=\pm 4\ \text{and\ } h_1\cdot H\equiv 0\mod 2.
\tag{3.1.36}
$$
In particular, if $(a,b)$ is a solution of $a^2-db^2=8$, the surface
$X$ has a nef element $h_1$ with square $4$ and the structure of a quartic,
if the $h_1$ is very ample.

The existence of an element $h_1\in N(X)$ satisfying \thetag{3.1.36}
(for one of signs $+$ or $-$) is equivalent to $Y\cong X$ for a general
$X$ with $\rho(X)=2$. In particular, $Y\cong X$ if $X$ has
a structure of quartic with the linear system $|h_1|$ of planes
of even degree with respect to the hyperplane $H$ of the intersection
of quadrics (i. e. $H\cdot h_1\equiv 0\mod 2$).
It is even sufficient to have an element
$\widetilde{h}_1\in W^{(-2)}(N(X))(h_1)$ with
$\widetilde{h}_1\cdot H\equiv 0\mod 2$.
\endproclaim

\demo{Proof} Let $(a,b)$ be a solution of $a^2-db^2=\pm 8$. It follows
that $a$, $b$ are odd. We have $(a-b\mu)(a+b\mu)\equiv 0\mod 8$.
It follows that either $a-b\mu\equiv 0\mod 4$
or $a+b\mu\equiv 0\mod 4$. Changing $b$ to $-b$ if necessary we
can assume that $a\equiv b\mu\mod 4$. Then $2a\equiv 2b\mu\mod 8$ and
$h_1\in N(X)$. We have $(h_1)^2=4(a^2-db^2)/8=\pm 8/2=\pm 4$.
Vice versa, if $h_1\in N(X)$ satisfies \thetag{3.1.36},
it can be written in the form \thetag{3.1.35}, where $(a,b)$ satisfies 
$a^2-db^2=\pm 8$.
\enddemo

Additionally applying Theorem 2.2.3, we get the following simple
sufficient condition when $Y\cong X$ which is valid for $X$ with any
$\rho (X)$:

\proclaim{Corollary 3.1.9}
Let $X$ be a K3 surface.
Assume that $X$ is an intersection of three quadrics (more generally,
$X$ has a primitive polarization $H$ of degree 8).
Let $Y$ be a K3 surface which is the double covering
of the net $\bp^2$ of the quadrics defining $X$ ramified along the 
curve of degenerate quadrics (more generally, $Y$ is the moduli space of 
sheaves on $X$ with the Mukai vector $v=(2,H,2)$). 

Then $Y\cong X$ if there exists $h_1\in N(X)$ such that
the primitive sublattice $[H,h_1]_{\text{\pr}}$ in $N(X)$
generated by $H$ and $h_1$ has odd determinant and
$$
(h_1)^2=\pm 4\ \text{and\ } h_1\cdot H\equiv 0\mod 2.
\tag{3.1.37}
$$

This condition is necessary to have $Y\cong X$
if either $\rho(X)=1$, or $\rho(X)=2$ and
$X$ is a general K3 surface with its Picard lattice
(i. e. the automorphism group of the transcendental periods
$(T(X),H^{2,0}(X))$ is $\pm 1$).
\endproclaim

\demo{Proof} The cases $\rho (X)\le 2$ had been considered.
We can assume that $\rho (X)>2$. Let $N=[H,h_1]_{\pr}$. All considerations
above for $N(X)$ of $\rk N(X)=2$ will be valid for $N$.
We can construct an associated with $h_1$ solution
$\widetilde{h}\in N$ with $\widetilde{h}^2=2$ such that $H$ and
$\widetilde{h}$ satisfy conditions of Theorem 2.2.3 for $N(X)$
replaced by $N$. It is easy to see that the conditions (b) and (c)
will be still satisfied if we extend $f$ in (b) $\pm$ identically on
the orthogonal complement $N^\perp_{N(X)}$. It finishes the proof.
\enddemo

It seems, many known examples of $Y\cong X$ (e. g. \cite{3}, \cite{15}) 
follow from the corollary.  

A similar (to Theorem 3.1.8) statement for $Y$ is simpler.

\proclaim{Theorem 3.1.10} Under notations of Proposition 3.1.2,
assume that
\newline
$d\equiv 1 \mod 8$. Then
elements
$$
h_1=(ah+b\delta)/2\in N(Y),
\tag{3.1.37}
$$
where $(a,b)$ is any integral solution of $a^2-db^2=\pm 8$
are all elements $h_1\in N(Y)$ with
$$
(h_1)^2=\pm 4.
\tag{3.1.38}
$$
In particular, if $(a,b)$ is a solution of $a^2-db^2=8$, the surface
$X$ has a nef element $h_1$ with square $4$, and a structure of quartic,
if the $h_1$ is very ample.

The equality $\gamma (H)=1$ (equivalently, $d=-\det N(X)$ is odd)
and the existence of an elements $h_1\in N(Y)$ satisfying
\thetag{3.1.38} are
equivalent to $Y\cong X$ for a general $Y$ with $\rho(Y)=2$.
In particular, $Y\cong X$, if $\gamma (H)=1$ and
$Y$ has a structure of quartic.
\endproclaim

\demo{Proof} It is similar and simpler.
\enddemo

We remark that if $N(X)$ of Theorem 3.1.8 has elements $h$
with $h^2=2$, then
$N(X)\cong N(Y)$ by Proposition 3.1.3. By Theorem 3.1.10,
then $Y\cong X$ if $N(X)$ has an element $h_1$ with $h_1^2=\pm 4$ 
(one does not need the congruence $h_1\cdot H\equiv 0\mod 2$).

\smallpagebreak

From our point of view, the statements 3.1.8 --- 3.1.10 are
very interesting also because they {\it give some geometric
meaning of  elements $h_1$ of the Picard lattice with the negative
square $(h_1)^2=-4$.} This is well-known only for elements
$\delta$ of the Picard lattice with
$\delta^2=-2$: then $\delta$ or $-\delta$ is effective.

\subhead
3.2. Geometry of $X\cong Y$, and calculation of the sets
$\Da_+$ and $\Da_-$
\endsubhead
Here we study geometry of surfaces $X$ and $Y$ with $\rho =2$ when
general $X\cong Y$. Equivalently, $N(X)\cong N(Y)\cong
N^8_d\cong N^2_d$ where $d\in \Da=\Da_+\cup \Da_-$. Moreover, we
give an algorithm to calculate the sets $\Da_+$ and $\Da_-$. Using
this algorithm, we find the first elements of these sets.

\smallpagebreak

It is more convenient to work with the surface $Y$ and its nef element
$h$ with $h^2=2$ which defines the structure of double plane, and
its  Picard lattice $N(Y)\cong N^2_d$. To define $X$, we should show a
polarization $H\in N(Y)$ with $H^2=8$ which defines the structure of
intersection of three quadrics. They should be related by the
correspondence described in Sect. 3.1. To find $H$, we show the
corresponding $h_1$ in $N(X)$ with $h_1^2=\pm 4$ such that $H$ is
the associated with $h_1$.

We use notations of Proposition 3.1.2. Thus,
$$
N(Y)=[h,\,\alpha,\,(h+\alpha)/2],\ h^2=2,\  \alpha^2=-2d,\ h\perp \alpha.
\tag{3.2.1}
$$
Here $d=-\det(N(Y))>0$ and $d\equiv 1\mod 8$, if $d\in \Da$.
We assume that $h$ is $nef$.

\proclaim{Proposition 3.2.1} Given $d\in \Da$,
the lattice $N(X)\cong N(Y)$ has non-zero elements with zero square,
if and only if $d$ is a square. It happens only for $d=1$ or $d=9$.

If $d=1$ we have $1\in \Da_+\cap \Da_-$, the set
$\Exc(Y)=\{\alpha\}$ (one should change $\alpha$ to $-\alpha$,
if necessary), and the K\"ahler cone $\overline{\Ka (Y)}$ is
generated by $h$ and by the
elliptic curve $E=(h-\alpha)/2$ (up to the linear equivalence).
This is the case (b) of Proposition 1.3.2.
The surface $Y$ has only one nef primitive element $H$ with
$H^2=8$. It defines $X$.
The element $H=(5h-3\alpha)/2=5E+\alpha$ which corresponds to
the case (c) of Proposition 1.3.1.
The $H$ is associated with the ample element
$h_1=(3h-\alpha)/2=3E+\alpha\in N(Y)$ with $(h_1)^2=4$. It is
also associated with the element $\widetilde{h}_1=(h-3\alpha)/2$
with $(\widetilde{h}_1)^2=-4$.

If $d=9$, $9\in \Da_- - \Da_+\cap \Da_-$.
For this case, $\Exc(Y)=\emptyset$, the K\"ahler cone is generated
by the elliptic curves $E_1=(3h-\alpha)/2$ and $E_2=(3h+\alpha)/2$,
the surface $Y$ has only one nef element with square 2 which is $h$,
and it has exactly two nef (and then ample) primitive elements
$H=(5h\pm \alpha)/2$ with $H^2=8$. They both define $X$ as intersection
of three quadrics and are associated with two elements
$h_1=(h\pm \alpha)/2$ with $h_1^2=-4$.
\endproclaim

\demo{Proof} Let $N(Y)$ has an element $c=(xh+y\alpha)/2$ with
$c^2=(x^2-dy^2)/2=0$ where $x,y$ are not both zero. Then $d$ is a
square in $\bq$, and then $d=(d_1)^2$ for $d_1\in \bn$. If $d\in \Da$,
then one of equations $(a-d_1b)(a+d_1b)=8$ or $(a-d_1b)(a+d_1b)=-8$
has an integral solution $(a,b)$. Simple calculations show that then
$d=d_1=1$, if the first equation has solution; and $d_1=1$ or $d_1=3$,
if the second equation has a solution. Respectively, $d=1$ or $d=9$.

Assume that  $d=1$. All elements $r=(ah+b\alpha)/2\in N(Y)$
with $r^2=-2$ correspond to
solutions of the equation $a^2-b^2=-4$ with $a\equiv b\mod 2$.
It follows that $(a,b)=(0,\pm 2)$, and the set $\Exc(Y)=\{\alpha\}$.
Thus, $nef$ elements $z=(ah+b\alpha)/2 \in N(Y)$ are defined by
the conditions: $z^2=(a^2-b^2)/2 \ge 0$, $a=h\cdot z > 0$
(these conditions give that $z\in \overline{V^+(Y)}$), and
$-b=z\cdot \alpha \ge 0$. It
follows that $z=qh+pE$ with $q\ge 0$ and
$p\ge 0$ where $E=(h-\alpha)/2$. Thus $h$ and $E$ generate the
$\overline{\Ka(Y)}$.
The $nef$ element $E$ has $E^2=0$ and is primitive. It follows
\cite{10} that $E$ is an elliptic curve on $X$ up to the linear
equivalence: $|E|$ is an elliptic pencil on $Y$.

The equation $a^2-b^2=8$ has the only solutions $(a,b)=(\pm 3, \pm 1)$.
The equation $a^2-b^2=-8$ has the only solutions $(a,b)=(\pm 1, \pm 3)$.
The corresponding $nef$ elements $h_1\in N(Y)$ with $(h_1)^2=4$
are the only $h_1=(3h-\alpha)/2$.
The corresponding associated $nef$ elements $H\in N(Y)$
are the only $H=(5h-\alpha)/2=5E+\alpha$. It follows that
$H$ is ample since its orthogonal complement is generated by
$\delta$ with $\delta^2=-8d=-8$.

Let $\widetilde{h}=(xh+y\alpha)/2\in N(Y)$ has $\widetilde{h}^2=2$.
Like above, we then get that $\widetilde{h}=\pm h$. Thus, $h$
is the only $nef$ element of $Y$ with square $2$.

Assume that $\widetilde{H}\in N(Y)$ is a primitive ample element with
$\widetilde{H}^2=8$. The pair $(Y,\widetilde{H})$ defines then a K3
surface $X$ with the polarization of degree $8$. Since $d\in \Da$,
we then know, by Theorem 3.1.7, that the corresponding surface
$Y$ with the polarization of degree two is $(Y,h)$
(since $Y$ has the only one nef element $h$
with $h^2=2$). Thus, $\widetilde{H}$ should be one of the
found above associated nef solutions, which are all equal to $H$.

Let $d=9$. Here we argue similarly to $d=1$. Like above
(considering the equation
$a^2-9b^2=-4$), we can see that $\Exc(Y)=\emptyset$. Then
$\overline{\Ka (Y)}=\overline{V^+(Y)}$ coincides with the set of all
$nef$ elements of $Y$ which are then $(xh+y\alpha)/2$ where
$x^2-9y^2\ge 0$ and $x>0$. It follows that $\overline{\Ka (Y)}$
is generated by the elliptic curves
$E_1$ and $E_2$. Like above (considering the equation $a^2-9b^2=4$),
we can see that $h$ is the only $nef$ element of $Y$ with $h^2=2$.
Like above, we can see that the equation $a^2-9b^2=-8$ has the only
solutions $(a,b)=(\pm 1,\, \pm 1)$. They give the only two associated
$nef$ elements $H=(5h\pm \alpha)/2$ with $H^2=8$. Like above, the
elements $H$ are the only primitive $nef$ elements
(they are then ample) of $Y$ with $H^2=8$.
The found above elliptic curves $E_1$ and $E_2$ are
the only elliptic curves (up to linear equivalence) on $Y$. We have
$H\cdot E_1\ge 3$ and $H\cdot E_2\ge 3$. By Proposition 1.3.1,
$X$ is then intersection of three quadrics.

It proves the statement.
\enddemo

\smallpagebreak

We now assume that $d\in \Da$ and $d>9$. Then $d$ is not a square,
$N(X)$ has no non-zero elements with square $0$, and $X$ is
intersection of three quadrics by Proposition 1.3.1.

We remind \cite{1} that the lattice $N(Y)\cong N^2_d$
(for $d\in \bn$, $d\equiv 1\mod 4$ and $d$ which is not a square)
with the fixed elements $h$ and $\alpha$
can be considered as an {\it order}
$N(Y)\subset \bk=\bq(\sqrt{d})$ of the real quadratic field
$\bk=\bq(\sqrt{d})$. One should correspond to
$z=(xh+y\alpha)/2\in N(Y)$, the element
$z=(x+y\sqrt{d})/2 \in \bk$. Thus, $1=h$ and $\sqrt{d}=\alpha$.
Then $N(Y)$ becomes a full
$\bz$-submodule and a subring of $\bk=\bq(\sqrt{d})$.
The $d$ is equal to the {\it discriminant} of the order $N(Y)$.
The order $N(Y)$ is the ring of all
integers in $\bk$ if $d$ is square-free.

The intersection pairing in $N(Y)$ then corresponds to
the norm $N$ in $\bk$:
$$
 z^2/2=N(z)=z\overline{z}=(x^2-dy^2)/4
\tag{3.2.2}
$$
where $\overline{z}=(x-y\sqrt{d})/2$.
An element $\epsilon \in N(Y)$ is called a {\it unit} if
$N(\epsilon )=\pm 1$. A unit $\epsilon$ is {\it positive}
(respectively {\it negative})
if $N(\epsilon )=1$ (respectively,
$N(\epsilon )=-1$). By \thetag{3.2.2}, positive units
are in one to one correspondence with $\widetilde{h}\in N(Y)$ with
$\widetilde{h}^2=2$, and negative units are in one
to one correspondence with roots $r\in N(Y)$ with $r^2=-2$.
Each unit $\epsilon$ defines an
automorphism $z\to \epsilon z$ of $N(Y)$ as a module,
which preserves the intersection
pairing (i.e. $\epsilon\in O(N(Y))$) and is proper (i.e.
$\epsilon \in SO(N(Y))$),
if $N(\epsilon)=1$. If $N(\epsilon)=-1$, then $\epsilon$
multiplies the intersection pairing by $-1$ (we shall call it then as
an {\it anti-automorphism} of the lattice $N(Y)$, $\epsilon\in AO(N(Y))$),
and it is non-proper, $\det(\epsilon )=-1$, i. e.
$\epsilon \in AO^-(N(Y))$. The group
$SO(N(Y))\cup AO^-(N(Y))$ is then identified with the group of units.
The group of units is ${\pm 1}\epsilon_0^{\bz}$
where
$$
\epsilon_0=(s+t\sqrt{d})/2=(sh+t\alpha)/2,\ s,\,t>0,
\tag{3.2.3}
$$
is the {\it fundamental unit.}  It is defined uniquely.
To get the full group $O(N(Y))\cup AO(N(Y))$ of automorphisms and
anti-automorphisms of the lattice $N(Y)$, one needs to add
the {\it involution $\phi_h\in O(N(Y))$ of
the double plane structure of $Y$:}
$\phi_h(h)=h$ and $\phi_h(\alpha)=-\alpha$. It preserves the
intersection pairing and is non-proper. We remind (e.g. see \cite{1})
that there exists a very effective algorithm for finding the fundamental
unit which uses continuous fractions.

\smallpagebreak

Depending on the norm of the fundamental unit $\epsilon_0$, one has
a different geometry of $Y$ and $X$.

Assume that $N(\epsilon_0)=-1$. Then $r=\epsilon_0$ has $r^2=-2$,
and (by elementary considerations)
$$
\Exc(Y)=\{r=(sh+t\alpha )/2,\,\widetilde{r}=\phi_h(r)=(sh-t\alpha )/2\}.
\tag{3.2.4}
$$
Thus, the $nef$ cone $\Ka(Y)$ of $Y$ is the set of
$z=(xh+y\alpha)/2\in N(Y)\otimes \br$ such that
$z\cdot h=x>0$, $z^2=(x^2-dy^2)/2 > 0$,
$r\cdot z=(xs-ytd)/2\ge 0$ and $\widetilde{r}\cdot z=(xs+ytd)/2\ge 0$.
The inequalities $x>0$ and $z^2=(x^2-dy^2)/2 > 0$ follow
from the last two (we assume that $z\not=0$). Thus,
$$
\Ka(Y)=\{z=(xh+y\alpha)/2\in N(Y)\otimes \br\ |\  xs\ge ytd\ge -xs\}.
\tag{3.2.5}
$$
Any automorphism of $Y$ should preserve the set $\Exc(Y)$.
It follows that the action of $Aut(Y)$ on $N(Y)$ is then generated
by the involution $\phi_h$, which generates the full automorphism group 
$Aut(Y)$ of a general $Y$ with $N(Y)$.
Thus, up to action of the group generated by
$\pm 1$, $W^{(-2)}(N(Y))$ and $Aut(Y)$, any $z\in N(Y)\otimes \br$
with $z^2>0$ has a unique representative $z_0$ in the
{\it fundamental domain}
$$
\Ka(Y)^+=\{(xh+y\alpha)/2\in N(Y)\otimes \br\ |\ xs\ge ytd\ge 0\}.
\tag{3.2.6}
$$
Similarly, any element $\omega \in N(Y)\otimes \br$ with $\omega^2<0$
has a unique representative $\omega_0$ such that the orthogonal
line to $\omega_0$ intersects the angle $\Ka(Y)^+$, and
$-\alpha\cdot \omega_0\ge 0$. Let $\omega=(xh+y\alpha)/2\in N(Y)$
where $\omega^2=(x^2-dy^2)/2<0$ and $x\ge 0$, $y>0$.
Then the orthogonal complement to $\omega$ is
$\br \omega^\perp$ where
$$
\omega^\perp=(dyh+x\alpha)/2,\ \text{and\ } (\omega^\perp)^2=-d\omega^2.
\tag{3.2.7}
$$
The above conditions mean that
$\omega^\perp=(dyh+x\alpha)/2 \in \Ka(Y)^+$, equivalently
$$
ys\ge xt\ge 0.
\tag{3.2.8}
$$
Elements $\omega_0\in N(Y)$ with $\omega_0^2<0$ satisfying
\thetag{3.2.8} are called {\it reduced}.
E.g. $r=\epsilon_0$ is reduced, and
$$
\epsilon_0^\perp=(dth+s\alpha)/2,\ \ (\epsilon_0^\perp)^2=2d.
\tag{3.2.9}
$$

We shall use the hyperbolic distance in the hyperbolic
space (it is a line for our case)
$$
\La((N(Y))=V^+(N(Y))/\br_{++}
\tag{3.2.10}
$$
where $\br_{++}$ is the set of positive real numbers. The distance
is given by
$$
\cosh(\rho(\br_{++}z_1, \br_{++}z_2))=
{z_1\cdot z_2\over \sqrt{z_1^2z_2^2}}
\tag{3.2.11}
$$
for $z_1,\,z_2\in V^+(N(Y))$.
The $\Ka(Y)^+/\br_{++}$ is a closed interval of the hyperbolic line
with terminals $\br_{++}h$ and $\br_{++}\epsilon_0^\perp$. It follows
that  points $\br_{++}z$ where $z=(xh+y\alpha)/2$ and $x^2-dy^2>0$,
of this interval are
characterized by the properties that
$$
x>0,\ \text{and\ }
1={h^2\over \sqrt{h^2h^2}}\le {h\cdot z\over \sqrt{2z^2}}\le
{h\cdot \epsilon_0^\perp\over \sqrt{h^2 (\epsilon_0^\perp)^2}}.
\tag{3.2.12}
$$
Equivalently, we have
$$
1\le {x\over \sqrt{x^2-dy^2}}\le {t\sqrt{d}\over 2}.
\tag{3.2.13}
$$
The interval has only a finite number of elements of
$N(Y)$ with the fixed square, i.e. when $x^2-dy^2>0$ is fixed.
Similar condition for $\omega=(xh+y\alpha)/2$ with
$\omega^2<0$ to be reduced is then
$$
{1\over \sqrt{d}}\le {y\over \sqrt{-(x^2-dy^2)}}\le {t\over 2}.
\tag{3.2.14}
$$

Assume that $\widetilde{h}\in \Ka(Y)$ is a nef element with
$\widetilde{h}^2=2$. Then there exists an involution $\phi$ of
$N(Y)$ such that $\phi(\widetilde{h})=\widetilde{h}$ and $\phi$ is $-1$
on the $\widetilde{h}^\perp$ which is generated by the element
$\widetilde{\alpha}$ with $\widetilde{\alpha}^2=-2d<-2$.
The automorphism  $\phi$ then preserves $\Ka(Y)$, the set $\Exc(Y)$, and
$\phi (\alpha)=\widetilde{\alpha}$. It follows that
$\phi=\phi_h$ and $\widetilde{h}=h$. Thus, the surface $Y$ has
the only one structure of double plane. It is defined by the $h$.

Assume that $H\in \Ka(Y)$ is a nef (and then ample) primitive element
with $H^2=8$. Assume that $H^\prime\in \Ka(Y)$ is another such an
element. By Proposition 3.1.1, there exists an automorphism
$\phi\in O(N(Y))$ such that $\phi(H)=H^\prime$. Then $\phi$ preserves
the fundamental chamber $\Ka(Y)$ for $W^{(-2)}(N(Y))$ and preserves
$Exc(Y)=\{r, \widetilde{r}\}$. It follows that either $\phi$ is identical or
$\phi=\phi_h$ is defined by the involution of the double plane of $Y$.
It follows that
$$
H=(xh+y\alpha)/2,\ \phi_h(H)=H^\prime=(xh-y\alpha)/2 
\tag{3.2.15}
$$
are the only primitive nef elements of $Y$ with square $8$.
They define equivalent structures $X\cong (Y,H)\cong (Y,H^\prime)$
of intersection of three quadrics. Since we assume that
$d\in \Da$, the corresponding double plane is $(Y,h)$ since
$h$ is the only nef element of $Y$ with $h^2=2$. By \thetag{3.2.13},
the element $H$ (and similarly $H^\prime$)
is characterized by the conditions (where $x,y$ are integers):
$$
\split
&H=(xh+y\alpha)/2,\ x^2-dy^2=16,\ x\equiv y\equiv 1\mod 2,\\
&x>0,\ y>0,\ 4<x< 2t\sqrt{d}.
\endsplit
\tag{3.2.16}
$$

If $z=(a+b\sqrt{d})/2$ gives a solution $(a,b)$ of \thetag{3.1.31},
then $\epsilon_0z$ gives a solution of \thetag{3.1.32},
since $N(\epsilon_0)=-1$. It follows that both equations
\thetag{3.1.31} and \thetag{3.1.32}
have solutions. Thus, $d\in \Da_+\cap \Da_-$. The opposite
statement is also valid. Really, assume that
$h_1=(ah+b\alpha)/2$ gives a solution of \thetag{3.1.31} and
$h_1^\prime=(a^\prime h+b^\prime \alpha)/2$ gives a solution of
\thetag{3.1.32}. Since $h_1^2=4$ and $(h_1^\prime)^2=-4$, like
for the propositions 3.1.1 and 3.1.2, one can prove that there exists
an anti-automorphism $\phi\in AO(N(Y))$ such that $\phi(h_1)=h_1^\prime$.
If $\det (\phi)=1$ the anti-automorphism
$\phi^\prime=\phi_h\phi$ has $\det(\phi^\prime)=-1$. It then corresponds
to a unit $\epsilon$ with $N(\epsilon)=-1$. It follows that
$N(\epsilon_0)=-1$ for the fundamental unit $\epsilon_0$. Thus,
conditions $N(\epsilon_0)=-1$ and $d\in \Da_+\cap \Da_-$
are equivalent, if $d\in \Da$.

By Theorem 3.1.7, $H$ is the associated solution \thetag{3.1.29} or
\thetag{3.1.30} with a solution of the equation \thetag{3.1.31} or
\thetag{3.1.32} respectively.

Assume that $H$ is associated with a solution
$h_1=(ah+b\alpha)/2$ of the equation \thetag{3.1.31}, thus 
$(h_1)^2=4$. Then $x=a^2-4$, $y=ab$. By \thetag{3.2.16},
we then have (changing signs of $a$, $b$ if necessary):
$a>0$, $b>0$ and $4<a^2-4\le 2t\sqrt{d}$. It follows that
$h_1\in \Ka(Y)^+$ is nef and then ample element with square $4$
since for elements $h_1$ with $(h_1)^2=4$
the condition \thetag{3.2.13} is
$$
2\sqrt{2}\le a \le t\sqrt{2d}.
\tag{3.2.17}
$$
Thus, if $H$ is associated with $h_1=(ah+b\alpha)/2$ where
$a>0$ and $b>0$ and
$(h_1)^2=4$, then $h_1$ is ample and
$$
8< a^2\le 2t\sqrt{d}+4.
\tag{3.2.18}
$$
The condition \thetag{3.2.18} is much stronger than \thetag{3.2.17}.

Now assume that $H$ is associated with $h_1=(ah+b\alpha)/2$ having
$h_1^2=-4$. Then $x=a^2+4$, $y=ab$. By \thetag{3.2.16},
we then have (changing signs of $a$ and $b$, if necessary)
$a>0$, $b>0$ and $4<a^2+4\le 2t\sqrt{d}$. Thus,
$$
0<a^2\le 2t\sqrt{d}-4.
\tag{3.2.19}
$$
By \thetag{3.2.14}, $h_1^\perp\in \Ka(Y)^+$ (or $h_1$ is reduced)
if and only if $2\sqrt{2}\le b\le \sqrt{2}t$. Since $a^2-db^2=-8$,
this is equivalent to
$$
0< a^2\le 2dt^2-8.
\tag{3.2.20}
$$
The condition \thetag{3.2.19} is much stronger than \thetag{3.2.20}.

Thus, if $d\in \Da$, then either equation \thetag{3.1.31} has
a solution $(a,b)$ satisfying \thetag{3.2.18}, or equation
\thetag{3.1.32} has a solution $(a,b)$ satisfying \thetag{3.2.19}.
This gives an effective algorithm to find out if $d\in \Da$.
Then $d\in\Da_+\cap \Da_-$.

If $d\in \Da$, then there exist $nef$ (and then ample) elements
$\widetilde{h}_1\in N(Y)$ with $(\widetilde{h}_1)^2=4$.
Let us find them. Like for $h^2=2$ or $H^2=8$ above,
there exist only two such elements $\widetilde{h}_1$.
One of them is $\widetilde{h}_1=(ph+q\alpha)/2$ where $p>0$, $q>0$,
and another one is $\phi_h(\widetilde{h}_1)=(ph-q\alpha)/2$.
If there exists a solution $(a,b)$ of \thetag{3.1.31},
satisfying \thetag{3.2.18}, then $(p,q)=(a,b)$ and
$\widetilde{h}_1=h_1=(ah+b\alpha)/2$, since the $h_1$ is $nef$.
If there exists a solution $(a,b)$, $a>0$, $b>0$ of \thetag{3.1.32},
satisfying \thetag{3.2.19}, then $\widetilde{h}_1=-
\epsilon_0(a-b\sqrt{d})/2$ has $(\widetilde{h}_1)^2=4$,
$\widetilde{h}_1 \in \Ka(Y)^+$, and the $\widetilde{h}_1$ is then
$nef$ and ample.

Finally we have proved

\proclaim{Theorem 3.2.2} Let $d\in \bn$, $d\equiv 1 \mod 8$ and $d>9$
is not a square. Further we follow conditions and notations
of Proposition 3.1.2. Let $\epsilon_0=(s+t\sqrt{d})/2$,where $s>0$, $t>0$,
be the fundamental unit of the order $N(Y)$ (where $1=h$ and
$\sqrt{d}=\alpha$).Then $d\in \Da_+\cap \Da_-$, if and only if the norm
$N(\epsilon_0)=-1$ and either the equation $a^2-db^2=8$ has a solution
$(a,b)$ satisfying \thetag{3.2.18} or the equation $a^2-db^2=-8$
has a solution $(a,b)$ satisfying \thetag{3.2.19}. This gives
an effective algorithm of $d\in \Da_+\cap \Da_-$
(see the Program in Sect. 5). For $d\in  \Da_+\cap \Da_-$,
let $h_1=(ah+b\alpha)/2$ corresponds to such a solution
$(a,b)$ with $a>0$ and $b>0$. We have $h_1^2=4$, for the first case,
and $h_1^2=-4$
for the second.

Further we assume that $d\in\Da_+\cap \Da_-$. Then the surface $Y$
has a unique $nef$ (and then ample) element of degree two; it is equal
to $h$ and defines the structure of a double plane for $Y$ and the
involution $\phi_h$ of the double plane such that
${\phi_h}^\ast(h)=h$, ${\phi_h}^\ast(\alpha)=-\alpha$.
The set $\Exc(Y)=\{r=\epsilon_0=(sh+t\alpha)/2,\,\widetilde{r}=
\phi_h^\ast (r)=
(sh-t\alpha)/2\}$. For a general $Y$ with $N(Y)$ the
involution $\phi_h$ generates $Aut(Y)$.

The surface $Y$ has exactly two nef (and then ample) elements
$\widetilde{h}_1$ and $\phi_h^\ast(\widetilde{h}_1)$ with
square $4$ which define two isomorphic structures of quartic on $Y$.
Here $\widetilde{h}_1=h_1$ if $h_1^2=4$, and
$\widetilde{h}_1=-\epsilon_0(a-b\sqrt{d})/2$
if $h_1^2=-4$ and $h_1=(ah+b\alpha)/2$.

The surface $Y$ has exactly two nef (and then ample) primitive elements
$H=(xh\pm y\alpha)/2$ with square $8$. They are associated with
the element $h_1$ above:
$$
(x,y)=(a^2-4,\,ab),\ \text{if\ }a^2-db^2=8,
$$
and
$$
(x,y)=(a^2+4,\,ab),\ \text{if\ }a^2-db^2=-8.
$$
These ample elements $H$ define isomorphic structures $X=(Y,H)$
of intersection of three quadrics on $Y$ such that the corresponding
double plane is isomorphic to $(Y,h)$.
\endproclaim

\smallpagebreak

Now we assume that $d\in \Da$ and $N(\epsilon_0)=1$.
Then all units of the order $N(Y)$ have the norm 1.
It follows that
$N(Y)$ has no roots $r$ with $r^2=-2$, and the $nef$ cone is
$\Ka (Y)=\overline{V^+(N(Y))}$. 
It follows that any $z=(xh+y\alpha)/2 \in N(Y)$
with $z^2>0$ and $x>0$ is $nef$ and ample.
All $nef$ elements $\widetilde{h}\in N(Y)$ with
$\widetilde{h}^2=2$ are in one to one correspondence with units
$\epsilon\in \epsilon_0^\bz$ (we call them $nef$ too).
For this correspondence, the nef element $\epsilon(h)$ with
square $2$ corresponds to the $nef$ unit $\epsilon$. Here we identify 
units with automorphisms of the order $N(Y)$.

Each nef element $\epsilon(h)$ (or nef element with square $2$)
defines a structure of
double plane for $Y$ and the involution $\phi_{\epsilon(h)}$ of $Y$
such that $\phi_{\epsilon(h)}(\epsilon(h))=\epsilon(h)$ and
$\phi_{\epsilon(h)}$
is $-1$ on the orthogonal complement to $\epsilon(h)$ in $N(Y)$.
We have
$\phi_{\epsilon_1(h)}\phi_{\epsilon_2(h)}={\epsilon_1}^2{\epsilon_2}^{-2}$.
For general $Y$ with $N(Y)$, the involutions $\phi_{\epsilon(h)}$ generate
the full automorphism group $Aut(Y)$ (we shall see that soon).

For the hyperbolic distance \thetag{3.2.11}, the nef element
$h^\prime=\epsilon_0(h)$ with square $2$
is the nearest to $h$ comparing with other nef elements with square $2$.
It follows that $\br_{++}h$ and $\br_{++}h^\prime$
are the terminals of the closed interval
$$
I=[\br_{++}h,\ \br_{++}h^\prime]\subset V^+(N(Y))/\br_{++}\ 
\text{where\ }h^\prime=\epsilon_0(h) 
\tag{3.2.21}
$$
of the hyperbolic line $V^+(N(Y))/\br_{++}$. For general $Y$,
the interval $I$ gives the fundamental domain for the action of
$Aut(Y)$ on the line, and the action of $Aut(Y)$ is generated
by reflections $\phi_h$ and $\phi_{h^\prime}$ in the terminals
$\br_{++}h$ and $\br_{++}h^\prime$ of $I$.
Really, otherwise, there exists an automorphism $\phi$ of $Y$ such that
$\phi^\ast(h)=h^\prime$ and $\phi^\ast(h^\prime)=h$. Then
$\phi$ is the symmetry in the center of the interval $I$, and
$\phi^\ast(h+h^\prime)=h+h^\prime$ (here $(h+h^\prime)^2>0$),
and $\phi$ is equal to $-1$ on the
orthogonal complement $\bz \beta$ to $h+h^\prime$. Thus, $\phi^\ast$ of
$N(Y)$ is generated by the reflection $s_\beta$ in the root
$\beta\in N(Y)$. Since $Y$ is general,
the action of $s_\beta$ should be $\pm 1$ on
the discriminant group of $N(Y)$ (to be continued as $\pm 1$ to the
transcendental lattice $T(Y)$). This is only possible if either
$\beta^2=-2$ or the orthogonal complement to $\beta$ in $N(Y)$ is
generated by  $nef$ element $e$ with
the square two. The first case is impossible since
$N(Y)$ has no elements with square $-2$. The second one is impossible
because $\br_{++}e=\br_{++}(h+h^\prime)$ is the center of the
interval $I$ which has smaller distance from the terminal $\br_{++}h$
than $\br_{++}h^\prime$. It also proves that the automorphism group of
$Y$ is generated by the involutions $\phi_{\epsilon(h)}$.

From this description of the action of $Aut(Y)$ on $N(Y)$,
we get that the $nef$ elements $h$ and $h^\prime=\epsilon_0(h)$ are not
conjugate by $Aut(Y)$, but their orbits give all $nef$ elements with
square $2$ of $Y$. There are two classes ${\epsilon_0}^{2\bz}(h)$
and ${\epsilon_0}^{2\bz}(h^\prime)={\epsilon_0}^{1+2\bz}(h)$ of $nef$
elements with square $2$ of $Y$, up to automorphisms of $Y$.

On the other hand, there exists an involution
$\xi=\epsilon_0\phi_h\in O^+(N(S))$ which acts as a symmetry in
the center of the interval $I$. Then $\xi(h)=h^\prime$,
$\xi(h^\prime)=h$. Clearly, the involutions $\phi_{\epsilon(h)}$ and
$\xi$ generate the automorphism group $O^+(N(Y))$ where $O^+(N(Y))$
is the subgroup of $O(N(Y))$ of index two, which preserves the
angle $V^+(N(Y))$. Really, this group is transitive
on the set of all nef elements in $N(Y)$ with square 2,
and its stabilizer subgroup has the order two.

Further we use the interval $I$ or the angle $\br_{++}I$ over $I$
instead of $\Ka(Y)^+$ for the previous case.

Up to the action of $\pm 1$ and $Aut(Y)$,
any $z\in N(Y)\otimes \br$ with $z^2>0$ has a unique representative
$z_0$ in $\br_{++}I$.
Like for \thetag{3.2.12}, $z=(xh+y\alpha)/2\in \br_{++}I$,
if and only if $z^2>0$, $x>0$ and
$$
1={h^2\over \sqrt{h^2h^2}}\le {h\cdot z\over \sqrt{2z^2}}\le
{h\cdot \epsilon_0(h)\over \sqrt{h^2 (\epsilon_0(h))^2}},
\tag{3.2.22}
$$
equivalently,
$$
1\le {x\over \sqrt{x^2-dy^2}}\le {s\over 2}.
\tag{3.2.23}
$$
Simple calculations show that $\br_{++}z$ belongs to the left half of
$I$ (containing $\br_{++}h$), if and only if
$$
1\le {x\over \sqrt{x^2-dy^2}}\le {\sqrt{s+2}\over 2}.
\tag{3.2.24}
$$
Similarly, an element $\omega \in N(Y)\otimes \br$ with $\omega^2<0$
has a unique representative $\omega_0$ such that the orthogonal
line to $\omega_0$ intersects the angle $\br_{++}I$, and
$-\alpha\cdot \omega_0\ge 0$. This $\omega_0$ is called {\it reduced.}
Like for \thetag{3.2.14}, we have
that $\omega=(xh+y\alpha)/2\in N(Y)\otimes \br$ where
$\omega^2=(x^2-dy^2)/2<0$ and $x\ge 0$, $y>0$,
is reduced (or $\omega^\perp=(dyh+x\alpha)/2 \in \br_{++}I$), if and
only if
$$
{1\over \sqrt{d}}\le {y\over \sqrt{-(x^2-dy^2)}}\le {s\over 2\sqrt{d}}.
\tag{3.2.25}
$$
Also $\omega^\perp=(dyh+x\alpha)/2$ belongs to the left half of $I$ if
and only if
$$
{1\over \sqrt{d}}\le {y\over \sqrt{-(x^2-dy^2)}}\le
{\sqrt{s+2}\over 2\sqrt{d}}.
\tag{3.2.26}
$$

By Proposition 3.1.1, the group $O^+(N(Y))$ is transitive on the set of
all primitive ample elements $H\in N(Y)$ with $H^2=8$.
It follows that there exist exactly two such elements $H$ and
$H^\prime=\xi(H)$ which give rays $\br_{++}H$ and
$\br_{++}H^\prime$ of the interval $I$.
If $H=H^\prime$, then all ample primitive elements of $Y$ with
square $8$ are $Aut(Y)$-equivalent and define equivalent
structures of double plane for $Y$. Moreover, by Theorem 3.1.6, any
structure of double plane of $Y$ can be obtained by this way. On
the other hand, $h$ and $h^\prime$ are not $Aut(Y)$-equivalent, and
define then different structures of double plane for $Y$. It proves
that $H^\prime$ is different from $H$. We further assume that $H$
defines the structure of double plane equivalent to $h$, and then
$H^\prime$ defines the structure of double plane equivalent to $h^\prime$.
Assume that $H=(xh+y\alpha)/2$. Since $\br_{++}H\in I$, then
$x>0$, $y>0$, and by \thetag{3.2.23}
$$
4<x<2s.
\tag{3.2.27}
$$

Let $d \in \Da_+$. Then $H=(xh+y\alpha)/2$ is associated
with a solution $(a,b)$ of $a^2-db^2=8$. Then
$h_1=(ah+b\alpha)/2\in N(Y)$ has ${h_1}^2=4$,
and $x=a^2-4$, $y=ab$. We can suppose that $a>0$, $b>0$.
By \thetag{3.2.27},
$$
8<a^2<2s+4
\tag{3.2.28}
$$
and
$$
2\sqrt{2}< a <\sqrt{2s+4}.
\tag{3.2.29}
$$
By \thetag{3.2.24}, this is equivalent to the condition that $\br_{++}h_1$
belongs to the left half of $I$. Like for
$nef$ elements with square $2$ or $8$ above, we can prove that
there exist exactly two $nef$ elements of $Y$ with square $4$
which belong to the interval $I$. They are conjugate by the involution
$\xi$ and give two non-isomorphic structures of quartic on $Y$ (for
a general $Y$ with $N(Y)$). The element $h_1$ belongs to
the left half of $I$. The element $\xi(h_1)$ belongs to the right
half of $I$ and gives another non-isomorphic
structure of quartic on $Y$. It gives the associated
element $\phi_{h^\prime}(H)$ with square $8$ in the image of $I$ by the
symmetry in its right terminal $\br_{++}h^\prime$.

Now assume that $d\in D_-$. Then $H=(xh+y\alpha)/2$ is associated
with a solution $(a,b)$ of $a^2-db^2=-8$,
$h_1=(ah+b\alpha)/2\in N(Y)$ has ${h_1}^2=-4$
and $x=a^2+4$, $y=ab$. We can suppose that $a>0$, $b>0$.
By \thetag{3.2.27},
$$
0<a^2<2s-4.
\tag{3.2.30}
$$
This is equivalent (use $a^2-db^2=-8$) to
$$
{1\over \sqrt{d}}\le {b\over \sqrt{8}}<{\sqrt{s+2}\over 2\sqrt{d}}
\tag{3.2.31}
$$
which (by \thetag{3.2.26}) means that $\br_{++}{h_1}^\perp$
belongs to the left half of $I$.
Like for $nef$ elements with square $2$ or $8$ above, we can prove that
there exist exactly two elements of $N(Y)$ with square $-4$ such
that the corresponding orthogonal elements  belong to
the interval $I$ (or there are exactly two reduced elements with
square $-4$). They are conjugate by the involution
$\xi$ and are not $Aut(Y)$-equivalent (for
a general $Y$ with $N(Y)$). The element ${h_1}^\perp$ belongs to
the left half of $I$. The element $\xi({h_1}^\perp)$
belongs to the right half of $I$.
The element $\xi(h_1)$ gives the associated
element $\phi_{h^\prime}(H)$ with square $8$ in the image of $I$ by the
symmetry in its right terminal $\br_{++}h^\prime$.

Thus, we have proved

\proclaim{Theorem 3.2.3} Let $d\in \bn$, $d\equiv 1 \mod 8$ and $d>9$
is not a square. Further we follow conditions and notations
of  Proposition 3.1.2. Let $\epsilon_0=(s+t\sqrt{d})/2$, $s>0$, $t>0$
be the fundamental unit of the order $N(Y)$ (where $1=h$ and
$\sqrt{d}=\alpha$). Then $d\in \Da_\pm-\Da_+\cap \Da_-$ if and
only if the norm $N(\epsilon_0)=1$ and the equation
$a^2-db^2=\pm 8$ has a solution
$(a,b)$ satisfying \thetag{3.2.28} for $a^2-db^2=8$, and satisfying
\thetag{3.2.30} for the equation $a^2-db^2=-8$. This gives
an effective algorithm of $d\in \Da_{\pm}-\Da_+\cap \Da_-$
(see the Program in Sect. 5). For $d\in \Da_{\pm}-\Da_+\cap \Da_-$,
we denote $h_1=(ah+b\alpha)/2\in N(Y)$ where $(a,b)$ is such
solution with $a>0$ and $b>0$. We have $h_1^2=\pm 4$.

Further we assume that $d\in \Da_{\pm}-\Da_+\cap \Da_-$.
The surface $Y$ has two nef (and then ample) elements of degree two,
$h$ and $h^\prime=\epsilon_0(h)$ , which
define two double plane structures of  $Y$ and the
involutions $\phi_h$ and $\phi_{h^\prime}$ of the double planes.
The set $\Exc(Y)=\emptyset$. For a general $Y$ with $N(Y)$ these two
double plane structures of $Y$ are not isomorphic (i.e. $h$ and $h^\prime$
are not $Aut(Y)$-equivalent), but any double plane
structure of $Y$ is isomorphic to one of these two;
the involutions $\phi_h$ and $\phi_{h^\prime}$ generate $Aut(Y)$.
There exists an involution $\xi\in O^+(N(Y))=\epsilon_0\phi_h^\ast$
such that $\xi(h)=h^\prime$; the involution $\xi$ and $Aut(Y)$ generate the 
group $O^+(N(Y))$.

The surface $Y$ has exactly two nef (and then ample) primitive
elements $H$ and $H^\prime$ with square $8$ such that
the rays $\br_{++}H$ and $\br_{++}H^\prime$ are in between the rays
$\br_{++}h$ and $\br_{++}h^\prime$. For general $Y$ with $N(Y)$,
the elements $H$ and $H^\prime$ give two non-isomorphic structures
of intersection of three quadrics of $Y$. Any structure of
intersection of three quadrics of $Y$ is isomorphic to one of them.
One of these elements, say
$H=(xh+b\alpha)/2$, is associated with the element $h_1$ above:
$$
(x,y)=(a^2-4,\,ab),\ \text{if\ }d\in \Da_+ - \Da_+\cap \Da_-,
\tag{3.2.32}
$$
and
$$
(x,y)=(a^2+4,\,ab),\ \text{if\ }d\in \Da_- - \Da_+\cap \Da_-.
\tag{3.2.33}
$$
The corresponding to $H$ double plane structure of $Y$ is isomorphic
to the one defined by $h$. The corresponding to $H^\prime$ double plane
structure of $Y$ is isomorphic to the one defined by $h^\prime$.

The ray $\br_{++}h_1$, if $d\in \Da_+-\Da_+\cap \Da_-$
(respectively, the ray
$\br_{++}{h_1}^\perp$, if $d\in \Da_- - \Da_+\cap \Da_-$)
belongs to the left half, containing $\br_{++}h$, of the interval
$I=[\br_{++}h,\br_{++}h^\prime]$ of the hyperbolic line
$V^+(N(Y))/\br_{++}$. For general $Y$ with $N(Y)$, the elements
$h_1$ and $\xi(h_1)$ are not $Aut(Y)$-equivalent;
any element $\widetilde{h}_1\in N(Y)$ with $\widetilde{h}_1^2=\pm 4$
is $Aut(Y)$-equivalent to $\pm h_1$ or $\pm \xi(h_1)$.
In particular, $Y$ has exactly two non-isomorphic structures of quartic
(defined by $h_1$ and $\xi(h_1)$) if $d\in \Da_+-\Da_+\cap \Da_-$.
\endproclaim

Now we consider conditions when equations \thetag{3.1.31} and
\thetag{3.1.32}, i.e. $a^2-db^2=\pm 8$, have solutions locally.
Equivalently, when these equations have solutions in the ring
$\bz_p$ of $p$-adic integers for any prime $p$ (obviously,
they have solutions in $\br$). We assume that
$d\equiv 1 \mod 8$. By Theorem 3.1.8, existence of
a solution is equivalent to existence of an element
$h_1\in N_d^2$ with $h_1^2=\pm 4$. Since the lattice
$N_d^2$ is even, $h_1$ is primitive in $N_d^2$.
Very often the genus of the lattice $N_d^2$ has only one class. Then
these local conditions are sufficient for existence of solutions of
the equations \thetag{3.1.31} and \thetag{3.1.32}. We have

\proclaim{Proposition 3.2.4} Assume that $d\in \bn$ and $d\equiv 1\mod 8$.

Then equation $a^2-db^2=\pm 8$ has a solution in $\bz_p$
(equivalently, the lattice $N_d^2\otimes \bz_p$ has an element $h_1$ with
$h_1^2=\pm 4$) for any prime $p$, if and only if
$$
\left(\,\pm 2\,\over \,p\,\right)=1
\tag{3.2.34}
$$
for any odd prime $p\vert d$.
We remind that
$$
\left(\,2\,\over \,p\,\right)=(-1)^{\omega(p)}
\text{ where } \omega(p)={p^2-1\over 8},
$$
and
$$
\left(-2\,\over \,p\,\right)=(-1)^{\omega(p)+\varepsilon(p)}
\text{ where } \varepsilon(p)={p-1\over 2}.
$$

In particular, if the genus of the lattice $N_d^2$ has only one
class, $cl(N_d^2)=1$, then the equation $a^2-db^2=\pm 8$
has an integral solution, if and only if \thetag{3.2.34} is valid for any odd
prime $p\vert d$.
\endproclaim

\demo{Proof} If $(a,b)$ is a solution of $a^2-db^2=\pm 8$,
then $(a/2)^2\equiv \pm 2\mod d$. It follows \thetag{3.2.34}.

Now we assume that $d\equiv 1\mod 8$ and \thetag{3.2.34} is valid.
We denote $S_p=N_d^2\otimes \bz_p$

Since $\det S=-d$, the lattice $S_p$ is
unimodular if $p\nmid d$.

An unimodular $p$-adic lattice, which is even for $p=2$,
is defined by its rank and determinant.
It follows that $S_2\cong U\otimes \bz_2$ (see \thetag{1.1.2}).
The lattice $U\otimes \bz_2$
has primitive elements with any even square. In particular, it has
primitive elements with the square $\pm 4$.
For odd prime $p\nmid d$ the lattice
$S_p\cong \langle 2\rangle\oplus \langle -2d \rangle\cong
\langle \pm 4 \rangle \oplus \langle \pm 4 (-d) \rangle$. Thus, $S_p$ has
elements with square $\pm 4$.

Assume that odd $p\vert d$. Then $S_p\cong \langle 2\rangle \oplus
\langle -2d \rangle \cong \langle
(\pm 2)\cdot 2 \rangle \oplus\langle -2d \rangle$
because $\pm 2\in (\bz_p^\ast)^2$
if \thetag{3.2.34} is valid. It follows that $S_p$ has
elements with square $\pm 4$.

This proves the statement.
\enddemo

Assume that $d\in  \bn$ and $d\equiv 1\mod 8$. In Sect. 5: Appendix,
we give Program for GP/PARI calculator (version 1.38) which uses
Proposition 3.2.4 and Theorems 3.2.2 and 3.2.3 to check if $d\in \Da_\pm$
where $d>9$. The program first checks that $d$ is not a square.
Then it checks that the local condition \thetag{3.2.34} is
satisfied for one of signs $\pm$. If that is true, it calculates
the fundamental unit $\epsilon_0=(s+t\sqrt{d})/2$ and its norm
$N_\epsilon=N(\epsilon_0)$. If $N_\epsilon=-1$, we apply
Theorem 3.2.2. We find all odd positive $(a,b)$ satisfying
$a^2-db^2=8$ and \thetag{3.2.18}, and $a^2-db^2=-8$ and \thetag{3.2.19}.
To make the algorithm more efficient, we also use that
$a^2\equiv \pm 8\mod d$. If such $(a,b)$ don't exist, then
$d\notin \Da$. If such $(a,b)$ exist, then $d\in \Da_+\cap\Da_-$, and we
get the element $h_1=(ah+b\alpha)/2$ with $h_1^2=\pm 4$.
If $N_\epsilon=1$, we similarly apply Theorem 3.2.3 to check that
$d\in \Da_\pm-\Da_+\cap\Da_-$ and to find $h_1=(ah+b\alpha)/2$
with  $h_1^2=\pm 4$, if $d\in \Da_\pm-\Da_+\cap\Da_-$.

If $d\in \Da_+\cap\Da_-$, we also find (according to Theorem 3.2.2)
classes $r=\epsilon_0(h)$ with $r^2=-2$, $\widetilde{h}_1$ with
$(\widetilde{h}_1)^2=4$ and $H=(xh+y\alpha)/2$ with $H^2=8$.

If $d\in \Da_\pm-\Da_+\cap\Da_-$, we also find (according to Theorem
3.2.3) classes $h^\prime=\epsilon_0(h)$ with $(h^\prime)^2=2$
and $H=(xh+b\alpha)/2$ with $H^2=8$.

The program additionally calculates the class number $cl=cl(d)$ of
the lattice $N_d^2$. If $cl=1$, then the local condition \thetag{3.2.34}
is sufficient for $d\in \Da_\pm$. This can be used to check that
calculations are correct.

Using this program, we get

\proclaim{Theorem 3.2.5} The first elements (up to 2009) of
$\Da=\Da_+\cup \Da_-$ are
$$
\split
&1(\pm),\, 9(-),\, 17(\pm),\, 33(-),\, 41(\pm),\, 57(-),\,73(\pm),\,
89(\pm),\, 97(\pm),\, 113(\pm),\, 129(-),\\
&137(\pm),\, 153(-),\,161(+),\, 177(-),\,193(\pm),\,
201(-),\,209(-),\,217(+),\, 233(\pm),\\
&241(\pm),\,249(-),\,281(\pm),\,297(-),\, 313(\pm),\,
329(+),\, 337(\pm),\, 353(\pm),\,369(-),\\
&393(-),\,409(\pm),\,417(-),\,433(\pm),\, 449(\pm),\,457(\pm),\,489(-),
\,497(+),\, 513(-),\\
&521(\pm),\, 537(-),\,553(+),\,561(-),\,
569(\pm),\, 593(\pm),\, 601(\pm),\, 617(\pm),\, 633(-),\\
&641(\pm),\,649(-),\, 657(-),\,673(\pm),\, 681(-),\,
713(+),\,721(+),\,737(-),\, 753(-),\\
&769,(\pm),\, 801(-),\, 809(\pm),\,833(+),\,849(-),\, 857(\pm),\,
873(-),\,881(\pm),\, 889(+),\\
&913(-),\, 921(-),\, 929(\pm),\, 937(\pm),\,
953(\pm),\, 969(-),\,977(\pm),\,1017(-),\,1033(\pm),\\
&1041(-),\,1049(\pm),\,1057(+),\,1081(+),\,
1097(\pm),\, 1121(-),\, 1137(-),\,
1153(\pm),\\
&1161(-),\,1169(+),\,1177(-),\, 1193(\pm),\,
1201(\pm),\, 1217(\pm),\,1233(-),\,1241(\pm),\\
&1249(\pm),\, 1273(-),\,1289(\pm),\,1321(\pm),\,
1329(-),\, 1337(+),\,1353(-),\,1361(\pm),\\
&1377(-),\,1401(-),\,1409(\pm),\,1433(\pm),\,1441(-),\,1457(+),\,
1473(-),\, 1481(\pm),\\
&1497(-),\,1513(+),\,1529(-),\,1553(\pm),\,1561(+),\,1569(-),\,1577(-),\,
1609(\pm),\\
&1633(+),\,1649(\pm),\,1657(\pm),\,1673(+),\,1689(-),\,1697(\pm),\,
1713(-),\,1721(\pm),\\
&1737(-),\,1753(\pm),\,1777(\pm),\,1793(-),\,1801(\pm),\,1809(-),\,
1817(+),\,1841(+),\\
&1857(-),\,1873(\pm),\,1881(-),\,1889(\pm),\,1913(\pm),\,1921(\pm),\,
1969(-),\,1977(-),\\
&1993(\pm),\,2009(+)
\endsplit
$$
where we mark $d\in \Da$ by $+$ (respectively $-$) if
$d\in \Da_+-\Da_+\cap\Da_-$ (respectively $d\in \Da_- - \Da_+\cap\Da_-$),
and by $\pm$, if $d\in \Da_+\cap\Da_-$.
\endproclaim

Calculations using the Program give for the first 10 non-square 
elements of $\Da$ the following (where $\omega=(1+\sqrt{d})/2$):

\smallpagebreak

\noindent
$d=17\in \Da_+\cap \Da_-$:
$cl(d)=1$, $\epsilon_0=3+2\omega$, $N(\epsilon_0)=-1$,
$h_1=(3h+\alpha)/2$, $r=(8h+2\alpha)/2$,
$\widetilde{h}_1=(5h+\alpha)/2$, $H=(13h+3\alpha)/2$;

\smallpagebreak

\noindent
$d=33\in \Da_- - \Da_+\cap\Da_-$:
$cl(d)=1$, $\epsilon_0=19 +8\omega$,
$N(\epsilon_0)=1$, $h_1=(5h+\alpha)/2$, $h^\prime=(46h+8\alpha)/2$,
$H=(29h+5\alpha)/2$;

\smallpagebreak

\noindent
$d=41 \in \Da_+\cap \Da_-$:
$cl(d)=1$, $\epsilon_0=27 +10\omega$, $N(\epsilon_0)=-1$,
$h_1=(7h+\alpha)/2$, $r=(64h+10\alpha)/2$, $\widetilde{h}_1=h_1$,
$H=(45h+7\alpha)/2$;

\smallpagebreak

\noindent
$d=57 \in \Da_- - \Da_+\cap\Da_-$:
$cl(d)=1$, $\epsilon_0=131 +40\omega$, $N(\epsilon_0)=1$,
$h_1=(7h+\alpha)/2$, $h^\prime=(302h+40\alpha)/2$,
$H=(53h+7\alpha)/2$;

\smallpagebreak

\noindent
$d=73 \in \Da_+\cap \Da_-$:
$cl(d)=1$,
$\epsilon_0=943 +250\omega$, $N(\epsilon_0)=-1$,
$h_1=(9h+\alpha)/2$, $r=(2136 h+250\alpha)/2$,
$\widetilde{h}_1=h_1$, $H=(77h+9\alpha)/2$;

\smallpagebreak

\noindent
$d=89 \in \Da_+\cap \Da_-$:
$cl(d)=1$, $\epsilon_0=447 +106\omega$, $N(\epsilon_0)=-1$,
$h_1=(9h+\alpha)/2$, $r=(1000h+106\alpha)/2$,
$\widetilde{h}_1=(217h+23\alpha)/2$, $H=(85h+9\alpha)/2$;

\smallpagebreak

\noindent
$d=97 \in \Da_+\cap \Da_-$:
$cl(d)=1$, $\epsilon_0=5035 +1138\omega$,
$N(\epsilon_0)=-1$, $h_1=(69h+7\alpha)/2$, $r=(11208h+1138\alpha)/2$,
$\widetilde{h}_1=h_1$, $H=(4757h+483\alpha)/2$;

\smallpagebreak

\noindent
$d=113 \in \Da_+\cap \Da_-$:
$cl(d)=1$, $\epsilon_0=703 +146\omega$,
$N(\epsilon_0)=-1$, $h_1=(11h+\alpha)/2$, $r=(1552h+146\alpha)/2$,
$\widetilde{h}_1=h_1$, $H=(117h+11\alpha)/2$;

\smallpagebreak

\noindent
$d=129 \in \Da_- - \Da_+\cap\Da_-$:
$cl(d)=1$, $\epsilon_0=15371 +2968\omega$, $N(\epsilon_0)=-1$,
$h_1=(11h+\alpha)/2$, $h^\prime=(33710h+2968\alpha)/2$,
$H=(125h+11\alpha)/2$;

\smallpagebreak

\noindent
$d=137 \in \Da_+\cap \Da_-$:
$cl(d)=1$, $\epsilon_0=1595 +298\omega$, $N(\epsilon_0)=-1$,
$h_1=(35h+3\alpha)/2$, $r=(3488h+298\alpha)/2$,
$\widetilde{h}_1=(199h+17\alpha)/2$, $H=(1229h+105\alpha)/2$.

\smallpagebreak

\subhead
3.3. An application to moduli of $X$ and $Y$
\endsubhead
The results above can be interpreted from the point of view of
moduli of intersections of three quadrics in $\bp^5$.

By period map and local or global Torelli Theorem,
the moduli space of K3 surfaces $X$ which are intersections of three
quadrics in $\bp^5$ (more generally, K3 surfaces $X$ with a primitive
polarization $H$ of degree $H^2=8$)
are $19$-dimensional. Moreover, surfaces $X$ with $\rho(X)=\rho$ belong to
a $20-\rho$-dimensional submoduli space. If $X$ is general, then
$\rho(X)=1$, and the surface $Y$ cannot be isomorphic to $X$
because $N(X)=\bz H$ where $H^2=8$, and $N(X)$ does not have elements
with square $2$ which is necessary if $Y\cong X$. {\it Thus, if $Y\cong X$,
then $\rho (X)\ge 2$, and $X$ belongs to a codimension 1 submoduli
space of K3 surfaces which is a divisor in the moduli space
(up to codimension 2). The set $\Da$ labels connected components
of the divisor. Each $d\in \Da$, gives an irreducible and connected
codimension one moduli subspace of K3 surfaces $X$ with the Picard lattice
$N(X)=N^8_d\cong N^2_d$ (more generally, $N^8_d\subset N(X)$, but
the polarization $H\in N^8_d$, see Corollary 3.1.9); $Y\cong X$ for
all $X$ from this subspace.}
See \cite{7}, \cite{8} on corresponding results about
connected components of moduli of K3 surfaces
with condition on Picard lattice. See also \cite{2}.

E.g. it is well-known that $Y\cong X$ if $X$ has a line. This is
a divisorial condition on moduli of $X$. The corresponding
component is labeled by $d=17\in \Da$.
Really, let $l\in N(X)$ be the class of line. Then
the intersection matrix of $H$ and $l$ is
$$
\pmatrix
H^2&H\cdot l\\
l\cdot H& l^2
\endpmatrix =
\pmatrix
8&1\\
1&-2
\endpmatrix
$$
which has the determinant $-17$. Since $17\in \Da$, then $Y\cong X$.
The projection from the line $l$
gives an embedding of $X$ to ${\bp}^3$ as a quartic. The corresponding
hyperplane section is $h_1=H-l$, it has ${h_1}^2=4$. We have
$H\cdot h_1=7$ which is odd.
The reflection of $h_1$ in $l$ gives $\widetilde{h}_1=h_1+(h_1\cdot l)l=
h_1+3l=H+2l$. We have $(\widetilde{h}_1)^2=4$ and
$\widetilde{h}_1\cdot H=10\equiv 0\mod 2$. Then $Y\cong X$ by Corollary
3.1.9. Of course, there exists a classical direct geometric 
isomorphism between $X$ and $Y$ in this case.

Our results in Sect. 2 imply: {\it There exists an infinite number of
different divisorial conditions on moduli of intersections $X$ of three
quadrics in $\bp^5$ such that each of them implies $Y\cong X$.
All these divisorial conditions are labeled by elements of the
infinite set $\Da\subset \bn$.
The number $d=17\in \Da$ corresponds to the classical example
of $X$ containing a line.}

It seems, a direct geometric construction of the isomorphism
between $X$ and $Y$ is known only for the first $d=1$, $9$ and $17$.
For all other $d\in \Da$ we were managed to prove that $Y\cong X$
only using the fundamental Global Torelli Theorem for
K3 surfaces proved by I.I. Piatetskii-Shapiro and I.R. Shafarevich in
\cite{10}.

{\it It would be interesting to find all possible codimension two
(or bigger) conditions on moduli of intersections $X$ of three
quadrics in $\bp^5$ which imply $Y\cong X$ and which don't follow from
the divisorial conditions on moduli which were described above.}

\head
4. A general perspective
\endhead

Similar methods and calculations can be developed in the
following general situation.

Let $X$ and $Y$ are K3 surfaces,
$$
\phi:(T(X)\otimes \bq, H^{2,0}(X))\cong (T(Y)\otimes \bq, H^{2,0}(Y))
\tag{4.1}
$$
an isomorphism of their transcendental periods over $\bq$, and
$$
(a_1,H_1,b_1)^{\pm},\dots (a_k,H_k,b_k)^{\pm},
\tag{4.2}
$$
a sequence of types of isotropic Mukai vectors of sheaves
on K3, and $\pm$ shows the direction of the correspondence.

Similar methods and calculations can be applied to answer the
following question:

{\it When there exists a correspondence between $X$ and
$Y$ given by the sequence \thetag{4.2} of types of Mukai vectors,
which gives the isomorphism \thetag{4.1} between their transcendental
periods?}

In \cite{6} and \cite{9} sufficient and necessary conditions 
on \thetag{4.1} are given when there exists at least one such
a sequence \thetag{4.2} with coprime Mukai vectors $(a_i,H_i,b_i)$.

In this paper, we had considered the case when $Y=X$, $\phi=\pm \text{id}$,
the sequence \thetag{4.2} consists of one primitive Mukai vector
$(2,H,2)^+$ with $H^2=8$ where $+$ means
that $Y$ is the moduli space of sheaves on $X$ with the
Mukai vector $(2,H,2)$.

\head
5. Appendix: A program for GP/PARI calculator
\endhead

\noindent
$\backslash$$\backslash$for d which is N, d$\backslash$equiv 1 mod 8
and d$>$9 it finds out
\newline
$\backslash$$\backslash$if d$\backslash$in $\backslash$Da, and
it finds basic polarizations
\newline
default(compatible,3)
\newline
pprint("\$d=",d,"\$");$\backslash$
\newline
if(issquare(d)==1,pprint("d=",d," is square and it is not in
$\backslash$Da"),$\backslash$
\newline
u=factor(d);uu=u[,1]~;u=uu;kill(uu);k1=matsize(u)[2];$\backslash$
\newline
dplus=1;dminus=1;for(k=1,k1,if(kro(2,u[k])==$-$1,dplus=$-$1,);$\backslash$
\newline
if(kro($-$2,u[k])==$-$1,dminus=$-$1,));$\backslash$
\newline
if(dplus==$-$1,pprint("d=",d," is not in $\backslash$Da\_+
by local conditions"),);$\backslash$
\newline
if(dminus==$-$1,pprint("d=",d," is not in $\backslash$Da\_$-$
by local conditions"),);$\backslash$
\newline
if(dplus==$-$1\&dminus==$-$1,pprint("d=",d," is not in $\backslash$Da by
local conditions"),$\backslash$
\newline
cl=classno(d);pprint("cl(d)=",cl);$\backslash$
\newline
if(dplus==1,aa0=vector(d,k,0);f=1;$\backslash$
\newline
for(k=1,d$-$1,if(mod(k$^\wedge$2,d)==mod(8,d),aa0[f]=k;f=f+1,));$\backslash$
\newline
apl0=vector(f$-$1,k,aa0[k]);kill(aa0);kill(f),);$\backslash$
\newline
if(dminus==1,aa0=vector(d,k,0);f=1;$\backslash$
\newline
for(k=1,d$-$1,if(mod(k$^\wedge$2,d)==mod($-$8,d),aa0[f]=k;f=f+1,));$
\backslash$
\newline
amin0=vector(f$-$1,k,aa0[k]);kill(aa0);kill(f),);$\backslash$
\newline
eps=unit(d);neps=norm(eps);$\backslash$
\newline
s=2$\ast$real(eps)+imag(eps);t=imag(eps);$\backslash$
\newline
pprint("fund unit=",eps);$\backslash$
\newline
pprint("norm fund unit=",neps);$\backslash$
\newline
if(neps==$-$1,$\backslash$
\newline
gampl=0;gammin=0;f=matsize(apl0)[2];lim2=sqrt(4+2$\ast$t$\ast$sqrt(d));$
\backslash$
\newline
a1=$-$d;for(n=0,lim2/d+1,a1=a1+d;for(m=1,f,a=a1+apl0[m];$\backslash$
\newline
if(type(a/2)==1,,a2=(a$^\wedge$2$-$8)/d;$\backslash$
\newline
if(issquare(a2)!=1,,b=isqrt(a2);if(a$^\wedge$2$-$4$<$=2$\ast$t$\ast$sqrt(d),$
\backslash$
\newline
gampl=1;apl=a;bpl=b;pprint("\$h\_1=(",apl,"h+",bpl,"$
\backslash$alpha)/2\$"),)))));$\backslash$
\newline
if(gampl==0,f=matsize(amin0)[2];lim2=sqrt(2$\ast$t$\ast$sqrt(d)$-$4);$
\backslash$
\newline
a1=$-$d;for(n=0,lim2/d+1,a1=a1+d;for(m=1,f,a=a1+amin0[m];$\backslash$
\newline
if(type(a/2)==1,,a2=(a$^\wedge$2+8)/d;$\backslash$
\newline
if(issquare(a2)!=1,,b=isqrt(a2);if(a$^\wedge$2+4$<$=2$\ast$t$\ast$sqrt(d),$
\backslash$
\newline
gammin=1;amin=a;bmin=b;$\backslash$
\newline
pprint("\$h\_1=(",amin,"h+",bmin,"$\backslash$alpha)/2\$"),))))),);$
\backslash$
\newline
if(gampl==0\&gammin==0,pprint("d=",d," is not in $\backslash$Da"),$
\backslash$
\newline
pprint("d=",d," is in $\backslash$Da\_+$\backslash$cap $
\backslash$Da\_$-$");$\backslash$
\newline
pprint("\$r=(",s,"h+",t,"$\backslash$alpha)/2\$");$\backslash$
\newline
if(gampl==1,x=apl$^\wedge$2$-$4;y=apl$\ast$bpl;$\backslash$
\newline
pprint("\$$\backslash$widetilde{h}\_1=h\_1\$");$\backslash$
\newline
pprint("\$H=(",x,"h+",y,"$\backslash$alpha)/2\$"),$\backslash$
\newline
pp=$-$(s$\ast$amin$-$t$\ast$bmin$\ast$d)/2;qq=$-
$(t$\ast$amin$-$s$\ast$bmin)/2;$\backslash$
\newline
x=amin$^\wedge$2+4;y=amin$\ast$bmin;$\backslash$
\newline
pprint("\$$\backslash$widetilde{h}\_1=(",pp,"h+",qq,"$
\backslash$alpha)/2\$");$\backslash$
\newline
pprint("\$H=(",x,"h+",y,"$\backslash$alpha)/2\$"))),$\backslash$
\newline
gampl=0;gammin=0;$\backslash$
\newline
if(dplus==1,$\backslash$
\newline
f=matsize(apl0)[2];lim2=sqrt(2$\ast$s+4);$\backslash$
\newline
a1=$-$d;for(n=0,lim2/d+1,a1=a1+d;for(m=1,f,a=a1+apl0[m];$\backslash$
\newline
if(type(a/2)==1,,a2=(a$^\wedge$2$-$8)/d;$\backslash$
\newline
if(issquare(a2)!=1,,b=isqrt(a2);if(a$^\wedge$2$-$4$<$=s$\ast$2,$\backslash$
\newline
gampl=1;apl=a;bpl=b;$\backslash$
\newline
pprint("\$h\_1=(",apl,"h+",bpl,"$\backslash$alpha)/2\$"),))))),);$\backslash$
\newline
if(dminus==1,$\backslash$
\newline
f=matsize(amin0)[2];lim2=sqrt(2$\ast$s$-$4);$\backslash$
\newline
a1=$-$d;for(n=0,lim2/d+1,a1=a1+d;for(m=1,f,a=a1+amin0[m];$\backslash$
\newline
if(type(a/2)==1,,a2=(a$^\wedge$2+8)/d;$\backslash$
\newline
if(issquare(a2)!=1,,b=isqrt(a2);if(a$^\wedge$2+4$<$=s$\ast$2,$\backslash$
\newline
gammin=1;amin=a;bmin=b;$\backslash$
\newline
pprint("\$h\_1=(",amin,"h+",bmin,"$\backslash$alpha)/2\$"),))))),);$
\backslash$
\newline
if(gampl==0\&gammin==0,pprint("d=",d," is not in $\backslash$Da"),$
\backslash$
\newline
if(gampl==1,pprint("\$d=",d,"\$ is in \$$\backslash$Da\_+$-$$
\backslash$Da\_+$\backslash$cap$\backslash$Da\_$-$\$");$\backslash$
\newline
pprint("\$h$^\wedge$$\backslash$prime=(",s,"h+",t,"$
\backslash$alpha)/2\$");$\backslash$
\newline
x=apl$^\wedge$2$-$4;y=apl$\ast$bpl;$\backslash$
\newline
pprint("\$H=(",x,"h+",y,"$\backslash$alpha)/2\$"),$\backslash$
\newline
pprint("\$h$^\wedge$$\backslash$prime=(",s,"h+",t,"$
\backslash$alpha)/2\$");$\backslash$
\newline
pprint("\$d=",d,"\$ is in \$$\backslash$Da\_$-$ $-$ $
\backslash$Da\_+$\backslash$cap$\backslash$Da\_$-$\$");$\backslash$
\newline
x=amin$^\wedge$2+4;y=amin$\ast$bmin;pprint("\$H=(",x,"h+",y,"$
\backslash$alpha)/2\$"))))));
\newline

\enddocument

\end